\documentclass[]{dakota-article}
\pdfoutput=1
\usepackage{amsmath,amsfonts,graphics,amssymb,amsthm}
\usepackage[lined,ruled,linesnumbered]{algorithm2e}
\usepackage{mathtools}
\newcounter{rmrk}[section]

\newenvironment{remark}{\stepcounter{rmrk}
\noindent {\newline \bf Remark
\arabic{section}.\arabic{rmrk}}}{\newline}
\numberwithin{equation}{section}
\newtheorem{cor}{Corollary}[section]
\newtheorem{theorem}{Theorem}[section]

\newtheorem{proposition}{Proposition}[section]
\newtheorem{definition}{Definition}[section]

\newtheorem{assumption}{Assumption}

\numberwithin{equation}{section}

\newcommand{\abs}[1]{\left\vert#1\right\vert}
\newcommand{\set}[1]{\left\{#1\right\}}

\newcommand{\pspace}{\mathbf{\Lambda}}
\newcommand{\dspace}{\mathbf{\mathcal{D}}}
\newcommand{\pmeas}{\mu_{\pspace}}
\newcommand{\dmeas}{\mu_{\dspace}}
\newcommand{\pborel}{\mathcal{B}_{\pspace}}
\newcommand{\dborel}{\mathcal{B}_{\dspace}}

\DeclareMathOperator*{\supp}{supp}

\newcommand{\priormeas}{P_{\pspace}^{\text{prior}}}
\newcommand{\postmeas}{P_{\pspace}^{\text{post}}}
\newcommand{\priordens}{\pi_{\pspace}^{\text{prior}}}
\newcommand{\postdens}{\pi_{\pspace}^{\text{post}}}
\newcommand{\pfpriormeas}{P_{\dspace}^{Q(\text{prior})}}

\newcommand{\pfpriordens}{\pi_{\dspace}^{Q(\text{prior})}}
\newcommand{\pfpostdens}{\pi_{\dspace}^{Q(\text{post})}}
\newcommand{\obsmeas}{P_{\dspace}^{\text{obs}}}
\newcommand{\obsdens}{\pi_{\dspace}^{\text{obs}}}
\newcommand{\postdenssbayes}{\tilde{\pi}_{\pspace}^{\text{post}}}

\def\SNL{Optimization and Uncertainty Quantification Department, Sandia National Labs, Albuquerque, NM}
\def\UCD{Department of Mathematical and Statistical Sciences, University of Colorado Denver, Denver, CO 80202}

\def\mycorauthor{T. Wildey}
\corauthor{\mycorauthor}
\coremail{tmwilde@sandia.gov}
\def\mytitle{A Consistent Bayesian Formulation for Stochastic Inverse Problems Based on Push-forward Measures}
\title{\mytitle}
\date{\today}
\author{T.~Butler\thanks{\UCD}, J.~Jakeman\thanks{\SNL}, T.~Wildey\samethanks[2]}
\funding{Sandia National Laboratories is a multi-program laboratory managed and operated by Sandia Corporation, a wholly 
owned subsidiary of Lockheed Martin Corporation, for the U.S. Department of Energy's National Nuclear Security 
Administration under contract DE-AC04-94AL85000.}
\shortauthor{Butler, Jakeman, Wildey}
\shorttitle{Bayesian Push-forward based Inference}

\begin{document}
\pagestyle{myheader}
\maketitle


\begin{abstract}
We formulate, and present a numerical method for solving, an inverse problem for inferring parameters of a deterministic model from stochastic observational data (quantities of interest).
The solution, given as a probability measure, is derived using a Bayesian updating approach for measurable maps that finds a posterior probability measure, that when propagated through the deterministic model produces a  push-forward measure that exactly matches the observed probability measure on the data. Our approach for finding such posterior measures, which we call {\em consistent Bayesian inference}, is simple and only requires the computation of the push-forward probability measure induced by the combination of a prior probability measure and the deterministic model. We establish existence and uniqueness of observation-consistent posteriors and present stability and error analysis. We also discuss the relationships between consistent Bayesian inference, classical/statistical Bayesian inference, and a recently developed measure-theoretic approach for inference. Finally,  
analytical and numerical results are presented to highlight certain properties of the consistent Bayesian approach and the differences between this approach and the two aforementioned alternatives for inference.
\end{abstract}

\keywords{stochastic inverse problems, Bayesian inference, uncertainty quantification, density estimation}

\amsid{60H30, 60H35, 60B10}

\section{Introduction}\label{sec:intro}

Inferring information about the parameters of a model from observations 
is a common goal in scientific modelling and engineering design. 
Given a simulation model, one often seeks to determine the set of possible model parameter values 
which produce a set of model outputs that match some set of observed or desired values. 
This problem is typically ill-posed and consequently some form of regularization is often imposed
to obtain a reasonable solution. 
There are a number of variants of this conceptual inverse problem, 
each imposing different regularization and operating under varying assumptions.

Inverse problems are often formulated deterministically and are solved using optimization-based approaches. 
In this paper, we consider the problem of inferring probabilistic information on model parameters, $\lambda$, from
probabilistic information describing model outputs that are functions of these parameters, which we refer to as quantities of interest (QoI) and denote by $Q(\lambda)$. 
Specifically, given a probability density on $Q(\lambda)$, we seek to determine a probability measure, $P_\pspace$, on the space of model parameters, $\Lambda$,
that is {\em consistent} with both the model and the observed data in the sense 
that the push-forward measure of $P_\pspace$ determined by the computational model
matches the probability measure on the observations.

Bayesian methods~\cite{Bernardo1994,Robert2001,Gelman2013,Jaynes1998,Kaipio_S_ACIP_2007,MarzoukNajmRahn,Stuart_AN_2010} are by far the most popular means
of inferring probabilistic descriptions of model parameters from data. 
Bayes' theorem states that the distribution of model parameters conditioned on data $d$, 
known as the posterior $\pi(\lambda|d)$, is proportional to the product of the prior distribution 
on the model parameters $\pi(\lambda)$ and the assumed relationship between the data and the model 
parameters known as the likelihood $\pi(d|\lambda)$, i.e. $\pi(\lambda|d)\propto\pi(\lambda)\pi(d|\lambda)$. 
The Bayesian approach allows one to embed prior beliefs through the prior $\pi(\lambda)$ and 
quantifies the change in the deviation from this prior probability, induced by the data $d$.

An alternative approach to stochastic inference uses measure-theoretic principles to construct a probability measure 
on the model parameters that is consistent with the model and the data in the sense described above~\cite{BE1,BE2,BE3}.
This measure-theoretic formulation exploits the geometry of the QoI map to define a 
unique set-based solution to the inverse problem in the space of equivalence classes of solutions,
referred to as generalized contours in~\cite{BE1,BE2,BE3}. 

In this paper we propose another framework for stochastic inversion that combines various
aspects of both the Bayesian and measure-theoretic approaches. 
Our approach begins with a prior probability measure/density on the model parameters, which 
is then updated to construct a posterior probability measure/density that is consistent with 
both the model and the observational data in the sense described above.
Specifically, given a prior probability density on the parameters, $\priordens$,
and an observed probability density on the QoI, $\obsdens$, our posterior density takes the form
\begin{equation}\label{eq:postpdf}
\postdens(\lambda) = \priordens(\lambda)\frac{\obsdens(Q(\lambda))}{\pfpriordens(Q(\lambda))}, \quad \lambda \in \Lambda,
\end{equation}
which we define formally using the disintegration theorem in Section~\ref{sec:Bayes}.
Here, $\pfpriordens$ denotes the push-forward of the prior through the model and represents a forward propagation of uncertainty.
Each of the terms in \eqref{eq:postpdf} has a particular statistical interpretation:
\begin{itemize}
\item $\priordens$ is the same prior utilized in the statistical/classical Bayesian formulation and represents any prior information described as a relative likelihood of parameter values.
\item $\pfpriordens$ represents how the prior knowledge of relative likelihoods of parameter values defines a  relative likelihood of model outputs by treating $\lambda$ as the predictor variables and $Q(\lambda)$ as the response variables.
\item $\obsdens$ describes the relative likelihood that the output of the model corresponds to the observed data.
\end{itemize}
Computing the posterior density~\eqref{eq:postpdf} only requires the construction of the push-forward of the prior $\pfpriordens$. The prior and the observed densities are assumed {\em a priori}.
The existence and uniqueness of such a posterior requires two key assumptions on push-forward measures that are described in detail in Section~\ref{sec:solvability}.
Stability of the posterior with respect to perturbations in the observational data is straightforward to show and is discussed in Section~\ref{sec:Stability}.

While we leverage aspects of traditional Bayesian and measure-theoretic approaches, the framework, assumptions, and general approach proposed here is distinct in several ways.
Unlike a classical statistical Bayesian inference approach, which infers the posterior distribution using a stochastic map, e.g., $d=Q(\lambda)+\epsilon$,
where $\epsilon$ is an assumed probabilistic error model, our approach directly inverts the observed stochasticity of the data, described as a probability measure or density, using the deterministic map $Q(\lambda)$.
Another significant departure from the classical Bayesian approach is that our approach requires the propagation of the prior probability 
distribution through the model (i.e., computation of the push-forward measure of the prior), 
which is often referred to as a forward propagation of uncertainty.   
In this paper, we utilize ``default'' methods (Monte Carlo sampling, kernel density estimation) 
to propagate uncertainties through relatively simple models.
Although more efficient techniques for forward propagation of uncertainty, 
such as
global polynomial approximations computed using stochastic spectral methods \cite{GhanemSpanos, XiuKarniadakis, WanKarniadakis, MaitreGhanem1, GhanemRedhorse}, sparse grid collocation methods \cite{gerstner03,hegland2003,bungartz04,Ma:2009:AHS:1514432.1514547}, and Gaussian process models~\cite{rasmussen2006} or
multi-fidelity methods~\cite{Narayan_GX_SISC_2013,Ng_E_AIAA_2012} do exist,
we do not consider them here.
The scope of this work is focused on the introduction and mathematical justification of a new paradigm for stochastic inference and the development of simple computational algorithms to numerically approximate and generate samples from the consistent posterior.
We prove that this approach is stable to perturbations in either the prior or observed densities in Section~\ref{sec:Stability}.
In a future work, we will consider various computational improvements to the algorithm for sampling from the consistent solution that are also stable to perturbations in the specified densities. 

Our Bayesian approach is actually more similar to the measure-theoretic approach developed in~\cite{BE1,BE2,BE3}
since they both provide consistent solutions to the stochastic inverse problem.
However, there are significant differences between the two approaches both conceptually and computationally. 
For completeness, we provide a brief description of the measure-theoretic approach in Section~\ref{subsec:compmeas}
and describe some of the similarities and differences.
In summary, the two approaches make different assumptions and exploit different information to find a unique 
solution to the stochastic inverse problem.


The remainder of this paper is organized as follows. 
We begin, in Section~\ref{sec:problems}, by describing the stochastic inverse problem and
we provide the precise mathematical assumptions that we make
in order to solve the stochastic inverse problem.
We also provide a set-based derivation of the posterior and prove that
this posterior provides a consistent solution on certain $\sigma$-algebras.
In Section~\ref{sec:Bayes}, we derive a consistent Bayesian solution to the stochastic inverse problem
on the Borel $\sigma$-algebra using the disintegration theorem.
We also provide a figure for ease of reference on the various probability spaces and steps involved in constructing this unique solution.
In Section~\ref{sec:Stability}, we prove that this unique posterior is stable with respect to perturbations in either the prior or observed densities.
We discuss some computational considerations in Section~\ref{sec:comp} and we derive a relationship between the error in the approximation of the push-forward of the prior and the error in the posterior.
Numerical results are presented in Section~\ref{sec:numerics} to build intuition and to demonstrate our computational algorithms on a set of applications.
In Section~\ref{sec:compare}, we compare our approach with
the statistical Bayesian formulation and a measure-theoretic approach.  
A simple example is used to explore situations where the consistent and statistical Bayesian approaches give identical results and cases where they differ.
Our concluding remarks are given in Section~\ref{sec:conclusions}.

\section{A deterministic model and a stochastic inverse problem}\label{sec:problems}

Consider a deterministic model $M(Y,\lambda)$ giving a solution $Y(\lambda)$ that is an implicit function of parameters $\lambda\in\pspace \subset \mathbb{R}^n$. 
The set $\pspace$ represents the largest physically meaningful domain of parameter values, and, for simplicity, we assume that $\pspace$ is finite-dimensional.
In practice, modelers are often only concerned with computing  a relatively small set of quantities of interest (QoI), $\{Q_i(Y)\}_{i=1}^m$, where $Q_i\in\mathbb{R}$ is a functional dependent on the model solution $Y$.
Since $Y$ is a function of parameters $\lambda$, so are the QoI and in the following we make this dependence explicit by replacing $Q_i(Y)$ with $Q_i(\lambda)$. 

Given a set of QoI, we define the QoI map $Q(\lambda) := (Q_1(\lambda), \cdots, Q_m(\lambda))^\top:\pspace\to\dspace\subset\mathbb{R}^m$ where $\dspace  := Q(\pspace)$ denotes the range of the QoI map.  
In the following, we define an inverse problem, which uses observed data, which, in turn, can be compared to the simulated QoI in order to infer information on $\pspace$.

\subsection{Defining the stochastic inverse problem}
Assume $\pspace$ and $\dspace$ are finite-dimensional metric spaces
and define $(\pspace, \pborel, \pmeas)$ to be a measure space using the Borel $\sigma$-algebra $\pborel$ and measure $\pmeas$. 
Similarly let $(\dspace, \dborel, \dmeas)$ be a measure space using the Borel $\sigma$-algebra $\dborel$ and measure $\dmeas$ where $\dspace  := Q(\pspace)$ as defined above. 
We refer to $\pmeas$ and $\dmeas$ as volume measures.
The role of the volume measures is quite different than the role of a probability measure, but nonetheless proves vital in a variety of stochastic and non-stochastic analyses \cite{BE1, BE3, Billingsly, Ferenczi1997}.
In fact, when a probability measure is described in terms of a probability density function, a dominating measure (referred to here as a volume measure) is at least implicitly defined, e.g., see the Radon-Nikodym theorem in \cite{Billingsly}.
In finite-dimensional problems, most parametric probability measures (e.g., Gaussian, Beta, uniform, etc.) have densities written in a form that assumes the dominating measure is the standard Lebesgue measure on $\mathbb{R}^n$.
Changing the dominating measure for a given probability measure will change the form of the density function while leaving the probabilities of events from the $\sigma$-algebra unchanged. 


 
It is often the case that the QoI map $Q$ is at least piecewise smooth so that $Q$ defines a measurable mapping between the measurable spaces $(\pspace, \pborel)$ and $(\dspace, \dborel)$.
For any $A\in\dborel$, the measurability of $Q$ implies that 
\[Q^{-1}(A) = \left\{ \lambda \in \pspace \ | \ Q(\lambda) \in A \right\}\in\pborel.\]
Furthermore, given any $A\in\dborel$ and $B\in\pborel$, 
\begin{equation}\label{eq:mapprops}
  Q(Q^{-1}(A)) = A, \quad \text{and} \quad B \subseteq Q^{-1}(Q(B)), 
\end{equation}
although we note that in most cases $B\neq Q^{-1}(Q(B))$ even when $n=m$. 

Next, we assume that we are given an observed probability measure, $\obsmeas$, on $(\dspace,\dborel)$
that is absolutely continuous with respect to $\pmeas$ and therefore admits an observed probability density, $\obsdens$.
We seek a probability measure, $P_\pspace$ on $(\pspace,\pborel)$ that is absolutely continuous with respect to $\pmeas$ and admits a probability density, $\pi_\pspace$,
such that, when propagated through the model, agrees with the observed density on the data almost everywhere.
We use $P^{Q(P_\pspace)}_\dspace$ to denote the push-forward of $P_\pspace$ through $Q(\lambda)$ 
which satisfies 
\begin{equation*}
P^{Q(P_\pspace)}_\dspace(A) = P_\pspace(Q^{-1}(A)).
\end{equation*}
for all $A\in \dborel$,

We define the inverse problem formally as follows:
\begin{definition}[Consistency]\label{def:inverse-problem}
Given a probability measure $\obsmeas$ on $(\dspace, \dborel)$ that is absolutely continuous with respect to $\dmeas$ and admits a density $\obsdens$ ,
the stochastic inverse problem seeks a probability measure $P_\pspace$ on $(\pspace, \pborel)$ that is absolutely continuous with respect to $\pmeas$ and admits a probability density
$\pi_\pspace$,
such that the subsequent push-forward measure induced by the map, $Q(\lambda)$, satisfies
\begin{equation}\label{eq:invdefn}
P_\pspace(Q^{-1}(A)) = P^{Q(P_\pspace)}_\dspace(A) = \obsmeas(A),
\end{equation}
for any $A\in \dborel$.
\end{definition}

We refer to any probability measure that satisfies \eqref{eq:invdefn} as a {\em consistent} solution. 
In Section \ref{sec:solvability}, we provide the necessary and sufficient conditions under which a consistent solution of
the stochastic inverse problem exists. Of course, this solution may not be unique, i.e., there
may be multiple probability measures that are consistent (satisfy Definition~\ref{def:inverse-problem}). 
A unique solution may be obtained by imposing additional constraints or structure on the problem. 
In this paper, such structure is obtained by incorporating prior information to construct a unique
Bayesian solution to the stochastic inverse problem.
Following the Bayesian philosophy \cite{Tarantola}, we introduce a {\em prior} probability measure $\priormeas$ on $(\pspace,\pborel)$ that is absolutely continuous with respect to $\pmeas$ and admits a probability density $\priordens$.
The prior probability measure encapsulates the existing knowledge about the uncertain parameters.
The topic of choosing an appropriate prior probability measure is important, and while we describe the basic assumptions required of the prior for our purposes below, a full discussion of this topic is beyond the scope of this work.

Assuming that $Q$ is at least measurable, then $\priormeas$ and $Q$ together induce a push-forward measure $\pfpriormeas$ on $\dspace$, which is defined for all $A\in \dborel$,
\begin{equation}\label{eq:pfprior}
\pfpriormeas(A) = \priormeas(Q^{-1}(A)).
\end{equation}

\subsection{Solvability of the stochastic inverse problem}\label{sec:solvability}


We first observe that even if $\priormeas$ is described in terms of a density $\priordens$ with respect to $\pmeas$, it is not necessarily the case that $\pfpriormeas$ is absolutely continuous with respect to $\dmeas$.
In other words, $\pfpriordens$ may not exist with respect to $\dmeas$.
For example, if $\dmeas$ is taken to be the usual Lebesgue measure and the QoI map is constant, then $\dmeas(\mathcal{D})=0$ yet $\pfpriormeas(\mathcal{D})=1$.
Assumption \ref{assump:density} allows us to circumvent this technical difficulty.
\begin{assumption}\label{assump:density}
Either the measure $\dmeas$ on $(\dspace, \dborel)$ is defined as the push-forward volume measure,  or the push-forward volume measure is absolutely continuous with respect to $\dmeas$.
\end{assumption}

It is worth noting that in practice we prefer using QoI maps that exhibit sensitivities to parameters $\lambda$ since otherwise it is virtually impossible to make useful inferences in $\pspace$.
Specifically, if $m\leq n$ and the Jacobian of $Q$ is defined and full rank everywhere in $\Lambda$, then the push-forward volume measure $\dmeas$ is absolutely continuous with respect to the usual Lebesgue measure.
In such a case, we can take $\dmeas$ as the usual Lebesgue measure and define the push-forward density $\pfpriordens$ in the usual way.
In any case, the precise definition of $\dmeas$ is only required to prove the theoretical results in Section~\ref{sec:Bayes} and is never actually constructed or approximated in practice.

Assumption \ref{assump:density} provides the necessary and sufficient conditions under which any probability measure defined on $(\pspace, \pborel)$ that is absolutely continuous with respect to $\pmeas$ defines a push-forward probability measure on $(\dspace, \dborel)$ that is absolutely continuous with respect to $\dmeas$.
In other words, with this assumption, we have that for any $A\in\dborel$,
\begin{equation}\label{eq:pfprior2}
\pfpriormeas(A) = \int_A \pfpriordens d\dmeas = \int_{Q^{-1}(A)} \priordens d\pmeas = \priormeas(Q^{-1}(A)).
\end{equation}


The consistency requirement imposes a type of compatibility condition on the measures $\obsmeas$ for which solutions exist.
Specifically, $\obsmeas$ must be absolutely continuous with respect to the push-forward volume measure $\dmeas$.
This compatibility condition is enforced by the following assumption that relates the prior and observed measures and the model.
\begin{assumption}\label{assump:dom}
There exists a constant $C>0$ such that $\obsdens(q)\leq C\pfpriordens(q)$ for a.e.~$q\in\dspace$.
\end{assumption}

\begin{remark}
Assumption~\ref{assump:dom} immediately implies that the observed measure is absolutely continuous with respect to the push-forward measure defined by the prior.
The theoretical results in this paper only require this absolute continuity condition, however, constructing a numerical approximation of the posterior, e.g., generating samples from the posterior using standard rejection sampling, requires the stronger assumption.
\end{remark}

This is an assumption on the {\em computability} of $\obsmeas$ (and thus on the posterior as discussed in Section~\ref{sec:comp}) 
and on the predictive capabilities of the prior and the model.
This is not very restrictive and simply assumes that any event that we observe with non-zero probability will be predicted using the model and prior with non-zero probability.
This assumption is consistent with the convention in Bayesian methods to choose the prior to be as general as possible because if the prior indicates that the probability of an event is zero, then no amount of data indicating a non-zero probability of the event will ever be incorporated into the posterior measure.
Throughout the remainder of this paper we assume that Assumptions \ref{assump:density} and \ref{assump:dom} hold.

Another practical consideration is that $\dspace = Q(\pspace)$ is generally not a hyperrectangle ($m$-orthotope) in $\mathbb{R}^m$ while in many situations it is more natural to
measure observable data on a hyperrectangle $\hat{{\cal D}} \subset \mathbb{R}^m$.
In this case, we assume $\hat{{\cal D}} \subset {\cal D}$ which is an assumption on the model and $\pspace$.
In other words, we assume that we can define our model and $\pspace$ such that we can predict all of the observed data.
This assumption may fail in cases with significant model bias, but we leave this topic for future work.

\subsection{A set-based derivation of a posterior}\label{subsec:bayesbrief}

Recall from classical probability theory that if $P$ is a probability measure, $B$ is an event of interest, and $A$ is an arbitrary event, then
\[
	P(B) = P(B|A)P(A) + P(B|A^c)P(A).
\]
Furthermore, if $B\subseteq A$, then this reduces simply to
\[
	P(B) = P(B|A)P(A).
\]
We now identify these abstract events $B$ and $A$ in the context of our problem.
First, to avoid technical issues involved with defining conditional probabilities on sets of zero probability, we will restrict all consideration to events $B\in\pborel$ such that $\priormeas(B)>0$, which implies that $\pfpriormeas(Q(B))>0$.
We let $A:=Q^{-1}(Q(B))$. 
Then, since $B\subseteq Q^{-1}(Q(B))$, we see that the above formula can be used to compute the probability of $B$ if the probability of $A$ and the probability of $B$ conditioned on $A$ can be determined.

We are motivated to use $P(A)=\obsmeas(Q(A))$ since the definition of the push-forward measure implies that the consistency condition is met due to the fact that $Q(A) = Q(B)$. 
We must now determine $P(B|A)$.
Since the observed measure can offer no insight as to this value, we use the prior measure $\priormeas$, the classical Bayes' theorem for events, and again the fact that $B\subseteq A$, to obtain
\[
	\priormeas(B|A) = \frac{\priormeas(A|B)\priormeas(B)}{\priormeas(A)} = \frac{\priormeas(B)}{\priormeas(A)}. 
\]
By definition of $\pfpriormeas$, we see that $\priormeas(A)=\pfpriormeas(Q(B))$.
We therefore  arrive at a formal expression for the posterior,
\begin{equation}\label{eq:setbayes}
\postmeas(B) := \begin{cases}
	\priormeas(B) \frac{\obsmeas(Q(B))}{\pfpriormeas(Q(B))}, & \text{if } \priormeas(B)>0, \\
	0, & \text{otherwise}, 
	\end{cases}
\end{equation}
which is simply a form that exploits the classical Bayes' rule while incorporating the push-forward measure induced by the prior and the model.
It is easy to show that this form defines a probability measure on a specific measurable space.
\begin{proposition}
The posterior given by \eqref{eq:setbayes} defines a probability measure on $(\pspace, \mathcal{C}_\pspace)$ where
\[
	\mathcal{C}_\pspace := \set{Q^{-1}(A)\, : \, A\in\dborel}
\]
is the so-called {\em contour} $\sigma$-algebra on $\pspace$ induced by the map $Q$.
\end{proposition}
\begin{proof}
It is clear that the posterior is non-negative and the measure of the empty set is zero.  We note that for any $B \in \mathcal{C}_\pspace$ we have $B = Q^{-1}(Q(B))$ which implies $\postmeas(B) = \obsmeas(B)$ since $\priormeas(B) = \priormeas(Q^{-1}(Q(B))) = \pfpriormeas(Q(B))$.  Thus, countable additivity of the posterior on $(\pspace,\mathcal{C}_\pspace)$ follows from the countable additivity of $\obsmeas$ on $(\dspace,\dborel)$.  The fact that $\postmeas(\pspace) = 1$ follows immediately from $\pspace = Q^{-1}(\dspace)$ since we define $\dspace = Q(\pspace)$.
\end{proof}

In fact, by construction it follows immediately that \eqref{eq:setbayes} defines a consistent solution to the stochastic inverse problem if $\pborel$ is replaced by $\mathcal{C}_\pspace$.
Furthermore, since $Q$ is assumed measurable, we have that $\mathcal{C}_\pspace\subset\pborel$, and $\postmeas$ is a consistent solution whenever $\mathcal{C}_\pspace=\pborel$, e.g. as happens when $Q$ is a bijection.
We can also show that \eqref{eq:setbayes} defines a measure on a conditional $\sigma$-algebra.
\begin{proposition}
Given a fixed $A\in \dborel$, \eqref{eq:setbayes} defines a measure on $\left(Q^{-1}(Q(A)), \mathcal{C}_A\right)$ where
\[\mathcal{C}_A = \left\{ B\in \pborel \ | \ Q(B) = Q(A) \right\}.\]
\end{proposition}
\begin{proof}
Again, it is clear that the posterior is non-negative and we can easily add the empty set to $\mathcal{C}_A$ and see that the measure of the empty set is zero.
Let $\left\{B_i \right\}_{i=1}^{\infty}$ be a collection of pairwise disjoint sets in $\mathcal{C}_A$.
We note that $\mathcal{C}_A \subset \pborel$ which implies that the prior is countably additive on $(Q^{-1}(Q(A)), \mathcal{C}_A)$.
Since $Q(B_i) = Q(A)$ by definition of $\mathcal{C}_A$, we have $Q(\cup_{i=1}^\infty B_i) = Q(A)$ and
\begin{multline*}
\postmeas\left( \cup_{i=1}^\infty B_i \right) = \priormeas\left( \cup_{i=1}^\infty B_i \right) \frac{\obsmeas(Q(A))}{\pfpriormeas(Q(A))} \\ = \sum_{i=1}^\infty \left(\priormeas(B_i) \right)\frac{\obsmeas(Q(A))}{\pfpriormeas(Q(A))} = \sum_{i=1}^\infty \left( \postmeas(B_i) \right).
\end{multline*}
Thus, countable additivity of the posterior on $(Q^{-1}(Q(A)), \mathcal{C}_A)$ follows from the countable additivity of the prior on $(Q^{-1}(Q(A)), \mathcal{C}_A)$.
\end{proof}

While this form of the posterior has an intuitive derivation, it can fail to be a probability measure on $(\pspace, \pborel)$ since it may not be countably additive with respect to pairwise disjoint events taken from $\pborel$ that do not belong to either $\mathcal{C}_\pspace$ or $\mathcal{C}_A$.
However, the decomposition of $\pspace$ into the directions informed by the data and the directions orthogonal to the data which must be regularized by the prior, motivates the utilization of the disintegration
theorem in Section~\ref{sec:Bayes} which is used to construct a consistent posterior probability measure on $(\pspace,\pborel)$.

\section{A consistent Bayesian solution to the stochastic inverse problem}\label{sec:Bayes}
We begin with a technical, but necessary, result from measure theory known as the disintegration theorem \cite{Dellacherie_Meyer, BE3}, which we state below in the context of probability measures. 
\begin{theorem}\label{thm:disintegration}
Assume  $Q:\pspace\to\dspace$ is $\pborel$-measurable, $P_\pspace$ is a probability measure on $(\pspace, \pborel)$ and $P_\dspace$ is the pushforward measure of $P_\pspace$ on $(\dspace, \dborel)$.
There exists a $P_\dspace$-a.e. uniquely defined family of conditional probability measures $\set{P_q}_{q\in\dspace}$ on $(\pspace, \pborel)$ such that for any $A\in\pborel$, 
\begin{equation*}
P_q (A) =  P_q ( A\cap Q^{-1}(q) ) , 
\end{equation*}
so $P_q (\pspace\setminus Q^{-1}(q)) = 0$, and there exists the following disintegration of $P_\pspace$,
\begin{equation}\label{eq:disintegration1}
	P_\pspace(A)  = \int_{\dspace} P_{q}(A) \, dP_\dspace(q) = \int_{\dspace} \bigg(  \int_{A\cap Q^{-1}(q)} \, dP _{q} (\lambda) \bigg) \, dP_\dspace(q),
\end{equation}
for $A \in \mathcal{B}_{\mathbf{\Lambda}}$.
\end{theorem}

Alternatively, given a probability measure $P_\dspace$ on $(\dspace,\dborel)$, the specification of a family of conditional probability measures $\set{P_q}_{q\in\dspace}$ on $(\pspace,\pborel)$ (with the properties specified in the above form of the disintegration theorem) defines a probability measure $P_\pspace$ through \eqref{eq:disintegration1} whose pushforward measure is exactly given by $P_\dspace$. 
We use this to prove that when $\postmeas$ is interpreted as an iterated integral of the density given by \eqref{eq:postpdf}, then it is in fact a probability measure that solves the stochastic inverse problem.

First, we describe the structure of the conditional probability measures obtained from the prior measure by application of the disintegration theorem.
Applying the above disintegration theorem and Assumption \ref{assump:density} to the prior measure $\priormeas$, for any $A\in\pborel$, we have
\begin{equation}\label{eq:disintegrate_prior1}
\priormeas(A) = \int_{\dspace} \bigg(  \int_{A\cap Q^{-1}(q)} \, dP_q^{\text{prior}}(\lambda) \bigg) \, d\pfpriormeas(q),
\end{equation}
Bayes' theorem states that the conditional density is given by,
\[\pi(\lambda | q) = \pi(\lambda) \frac{\pi(q | \lambda)}{\pi(q)}.\]
The product of $\pi(\lambda)$ and $\pi(q|\lambda)$ defines a joint density, which we set to be $\priordens(\lambda)$.
Since Bayes' theorem is only being applied on $\lambda\in Q^{-1}(q)$, we see that $\pi(q)$ is given by
%
\[\pi(q) = \int_{\lambda \in Q^{-1}(q)} \priordens(\lambda)\ d\mu_{\pspace,q}(\lambda) = \pfpriordens(q),\]
where $d\mu_{\pspace,q}$ is defined as the disintegration of the volume measure $\pmeas$ (for details on disintegrations of volume and other measures  see \cite{BE3}).
Since $\lambda\in Q^{-1}(q)$, we have that $q=Q(\lambda)$, and we can make the appropriate substitution to obtain
\[
\pi(\lambda | q) = \frac{\priordens(\lambda)}{\pfpriordens(Q(\lambda))} \quad \text{ and } \quad
dP_q^{\text{prior}}(\lambda) = \frac{\priordens(\lambda)}{\pfpriordens(Q(\lambda))} d\mu_{\pspace,q}(\lambda).
\]
Then, with Assumption \ref{assump:dom} and letting $dP_\dspace(q) = \obsdens(q)\, d\dmeas(q)$, we can use the conditionals defined above along with the iterated integral \eqref{eq:disintegration1} from the disintegration theorem to define a measure $P_\pspace$ that is a consistent solution to the stochastic inverse problem as  
\begin{equation}\label{eq:consistent_meas}
	P_\Lambda(A) = \int_{\dspace} \bigg(  \int_{A\cap Q^{-1}(q)} \, \frac{\priordens(\lambda)}{\pfpriordens(Q(\lambda))} d\mu_{\pspace,q}(\lambda) \bigg) \obsdens(q)\, d\dmeas(q). 
\end{equation}
This can be rewritten as
\begin{equation}\label{eq:consistent_meas1}
	P_\Lambda(A) = \int_{\dspace} \bigg(  \int_{A\cap Q^{-1}(q)} \, \priordens(\lambda)\frac{\obsdens(q)}{\pfpriordens(Q(\lambda))} d\mu_{\pspace,q}(\lambda) \bigg) \, d\dmeas(q). 
\end{equation}
Finally, for all $\lambda\in \set{A\cap Q^{-1}(q)}$, $Q(\lambda)=q$, and we can rewrite this as
\begin{equation}\label{eq:consistent_meas2}
	P_\Lambda(A) = \int_{\dspace} \bigg(  \int_{A\cap Q^{-1}(q)} \, \priordens(\lambda)\frac{\obsdens(Q(\lambda))}{\pfpriordens(Q(\lambda))} d\mu_{\pspace,q}(\lambda) \bigg) \, d\dmeas(q). 
\end{equation}
In this form, it is now evident that $P_\Lambda$ is formally the posterior $\postmeas$ interpreted in terms of the above iterated integral of the density given in \eqref{eq:postpdf}.
In other words, proper interpretation of the posterior probability measure as an iterated integral defines a solution to the stochastic inverse problem that is consistent in the sense of \eqref{eq:invdefn}.
This proves the following
\begin{theorem}\label{thm:consistent_posterior}
The probability measure $\postmeas$ on $(\pspace, \pborel)$ defined by
\begin{equation}\label{eq:consistent_posterior}
	\postmeas(A) = \int_{\dspace} \bigg(  \int_{A\cap Q^{-1}(q)} \, \priordens(\lambda)\frac{\obsdens(Q(\lambda))}{\pfpriordens(Q(\lambda))} d\mu_{\pspace,q}(\lambda) \bigg) \, d\dmeas(q), \ \forall A\in\pborel
\end{equation}
is a consistent solution to the stochastic inverse problem in the sense of \eqref{eq:invdefn}.
\end{theorem}

\begin{figure}[htbp]
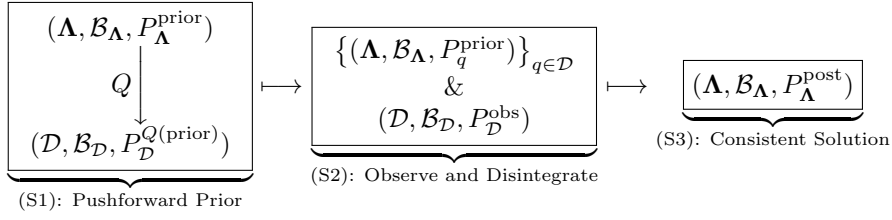

\begin{equation*}
\underbrace{\boxed{\begin{array}{c}
	(\pspace, \pborel, \priormeas) \\
	Q\Bigg\downarrow \\
	(\dspace,\dborel,\pfpriormeas)
\end{array}}}_{\text{(S1): Pushforward Prior}}
\longmapsto
\underbrace{\boxed{\begin{array}{c}
\set{(\pspace, \pborel, P_{q}^{\text{prior}})}_{q\in\dspace} \\
	\& \\
	(\dspace,\dborel,\obsmeas)
\end{array}}}_{\text{(S2): Observe and Disintegrate}}
\longmapsto
\underbrace{\boxed{(\pspace,\pborel,\postmeas)}}_{\text{(S3): Consistent Solution}}
\end{equation*}
\caption{The first step (S1) involves using the prior measure and QoI map to construct the pushforward of the prior on $(\dspace,\dborel)$. The second step (S2) involves constructing the conditional probability measures from the prior measure and its pushforward and determining the observed density. Finally, we use this information to construct a consistent solution in the third step (S3). }
\label{fig:process}
\end{figure}

For ease of reference, in Figure~\ref{fig:process}, we summarize the probability spaces involved in defining the consistent solution given by Theorem~\ref{thm:consistent_posterior}.
First, we use the map $Q$ to construct the pushforward of the prior measure.
Then, we construct a family of conditional probability densities that can be combined with an observed probability measure on the outputs of $Q$.
Finally, by the disintegration theorem, the combination of the observed probability measure and conditional probabilities forms a consistent solution.  

We note that the posterior density does not contain a normalization constant.
The fact that the posterior density integrates to one (equivalently that $\postmeas(\pspace)=1$)
follows from the utilization of the disintegration theorem, but for the sake of clarity we include the following result.

\begin{cor}
The posterior measure of $\pspace$ is 1.
\end{cor}
\begin{proof}
Using the definition of the posterior measure given by~\eqref{eq:consistent_posterior}, we have
\begin{eqnarray*}
\postmeas(\pspace) &=& \int_{\dspace} \bigg(  \int_{\pspace \cap Q^{-1}(q)} \, \priordens(\lambda)\frac{\obsdens(q)}{\pfpriordens(q)} d\mu_{\pspace,q}(\lambda) \bigg) \, d\dmeas(q) \\
&=& \int_{\dspace} \frac{\obsdens(q)}{\pfpriordens(q)} \bigg(  \int_{\pspace \cap Q^{-1}(q)} \, \priordens(\lambda) d\mu_{\pspace,q}(\lambda) \bigg) \, d\dmeas(q) \\
&=& \int_{\dspace} \obsdens(q) \, d\dmeas(q) \\
&=& 1
\end{eqnarray*}
where we use the fact that $Q(\lambda)=q$ and is constant on $\pspace \cap Q^{-1}(q)$, the definition of the push-forward density, and the fact that the observed density integrates to 1.
\end{proof}


\section{Stability of the consistent solution}\label{sec:Stability}

We have so far only considered the existence and uniqueness of the consistent solution to the stochastic inverse problem. 
Here, we consider the {\em stability} of this solution to perturbations in the specified observed and prior densities.
First, we require a notion of distance for probability measures $P$ and $Q$ defined on $(\pspace,\pborel)$.
Assuming that $P$ and $Q$ are absolutely continuous with respect to $\pmeas$, with densities given by $\pi_P$ and $\pi_Q$, respectively, we use the total variation distance (the so-called ``statistical distance'' metric), to compute
\begin{equation}\label{eq:TV_metric}
	d_{TV}(P,Q) := \int_\pspace \abs{\pi_P-\pi_Q}\, d\pmeas.
\end{equation}
The total variation metric is common in probability theory, but there are other choices we could use, e.g., see \cite{ChoosingMetrics_Gibbs}.
However, the total variation metric is one of the more straightforward metrics to compute, e.g., using Monte Carlo techniques. 
\begin{definition}
Given $\priormeas$ and $\obsmeas$, let $\widehat{\obsmeas}$ be any perturbation to $\obsmeas$ on $(\dspace,\dborel)$ satisfying Assumption~\ref{assump:dom}.
Let $\postmeas$ and $\widehat{\postmeas}$ denote the consistent solutions associated with $\obsmeas$ and $\widehat{\obsmeas}$, respectively.
We say that $\postmeas$ is stable with respect to perturbations in $\obsmeas$ if for all $\epsilon>0$ there exists $\delta>0$ such that
\begin{equation}\label{eq:stability_obs}
	d_{TV}(\obsmeas,\widehat{\obsmeas}) <\delta \Rightarrow     d_{TV}(\postmeas,\widehat{\postmeas}) < \epsilon.
\end{equation}
\end{definition}
\begin{theorem}\label{thm:stability_obs}
The consistent solution to the stochastic inverse problem $\postmeas$ is stable with respect to perturbations in $\obsmeas$.
\end{theorem}
\begin{proof}
Let $\epsilon>0$ be given and suppose $\widehat{\obsmeas}$ is chosen so that
$$
	d_{TV}(\obsmeas,\widehat{\obsmeas})<\epsilon.
$$
Then,
\begin{eqnarray*}
	d_{TV}(\postmeas,\widehat{\postmeas}) &=& \int_{\dspace} \bigg(  \int_{\pspace\cap Q^{-1}(q)} \, \frac{\priordens(\lambda)}{\pfpriordens(Q(\lambda))}\abs{\obsdens(Q(\lambda))-\widehat{\obsdens}(Q(\lambda))} d\mu_{\pspace,q}(\lambda) \bigg) \, d\dmeas(q) \\
	 &=& \int_{\dspace} \bigg(  \int_{\pspace\cap Q^{-1}(q)} \, \frac{\priordens(\lambda)}{\pfpriordens(Q(\lambda))} d\mu_{\pspace,q}(\lambda) \bigg) \abs{\obsdens(q)-\widehat{\obsdens}(q)}\, d\dmeas(q).
\end{eqnarray*}
The inner integral is taken over $\pspace\cap Q^{-1}(q)$, which represents the entire support of the conditional density, so is equal to $1$ for each $q$.
Thus,
\begin{equation*}
	d_{TV}(\postmeas,\widehat{\postmeas}) = \int_\dspace \abs{\obsdens(q)-\widehat{\obsdens}(q)}\, d\dmeas(q).
\end{equation*}
Since the integral on the right is equal to $d_{TV}(\obsmeas,\widehat{\obsmeas})$, the conclusion follows. 
\end{proof}

It is evident from Eq.~\eqref{eq:disintegration1} that perturbations to the prior measure at values of $\lambda$ where $\obsdens(Q(\lambda))=0$ do not impact the solution to the stochastic inverse problem. 
This observation motivates the following definition of stability with respect to perturbations in $\priormeas$.
\begin{definition}
Given $\priormeas$ and $\obsmeas$, let $\widehat{\priormeas}$ be any perturbation to $\priormeas$ such that Assumption~\ref{assump:dom} is still valid.
Let $\postmeas$ and $\widehat{\postmeas}$ denote the consistent solutions associated with $\priormeas$ and $\widehat{\priormeas}$, respectively.
Let $\set{P_q^{\text{prior}}}_{q\in\dspace}$ and $\set{\widehat{P_q^{\text{prior}}}}_{q\in\dspace}$ be the conditional probabilities defined by the disintegrations of $\priormeas$ and $\widehat{\priormeas}$, respectively. 
We say that $\postmeas$ is stable with respect to perturbations in $\priormeas$ if for all $\epsilon>0$ there exists $\delta>0$ such that for almost every $q\in\supp \obsdens$, 
\begin{equation}\label{eq:stability_prior}
	d_{TV}(P_q^{\text{prior}},\widehat{P_q^{\text{prior}}}) <\delta \Rightarrow     d_{TV}(\postmeas,\widehat{\postmeas}) < \epsilon.
\end{equation}
\end{definition}
\begin{theorem}\label{thm:stability_prior}
The consistent solution to the stochastic inverse problem $\postmeas$ is stable with respect to perturbations in $\priormeas$.
\end{theorem}
\begin{proof}
Let $\epsilon>0$ be given and suppose $\widehat{\priormeas}$ is chosen such that for a.e.~$q\in\supp \obsdens$, 
$$
	d_{TV}(P_q^{\text{prior}},\widehat{P_q^{\text{prior}}}) < \epsilon.
$$
Then, using the more compact notation of $\pi_q^{\text{prior}}$ and $\widehat{\pi_q^{\text{prior}}}$ to represent the conditional densities associated with $P_q^{\text{prior}}$ and $\widehat{P_q^{\text{prior}}}$, respectively, we have
\begin{eqnarray*}
	d_{TV}(\postmeas,\widehat{\postmeas}) &=& \int_{\dspace} \bigg(  \int_{\pspace\cap Q^{-1}(q)} \, \abs{\pi_q^{\text{prior}} - \widehat{\pi_q^{\text{prior}}}} \obsdens(Q(\lambda)) d\mu_{\pspace,q}(\lambda) \bigg) \, d\dmeas(q) \\
	 &=& \int_{\dspace} \bigg(  \int_{\pspace\cap Q^{-1}(q)} \, \abs{\pi_q^{\text{prior}} - \widehat{\pi_q^{\text{prior}}}} d\mu_{\pspace,q}(\lambda) \bigg) \obsdens(q)\, d\dmeas(q) \\
	 &=& \int_{\dspace} d_{TV}(P_q^{\text{prior}},\widehat{P_q^{\text{prior}}}) \obsdens(q)\, d\dmeas(q) \\
	 &=&\int_{q\in \supp \obsdens} d_{TV}(P_q^{\text{prior}},\widehat{P_q^{\text{prior}}}) \obsdens(q)\, d\dmeas(q) \\
	 &<& \epsilon \int_{q\in \supp \obsdens} \obsdens(q)\, d\dmeas(q) \\
	 &=& \epsilon. 
\end{eqnarray*}
\end{proof}

We conclude with some brief interpretations of these stability results.
Theorem~\ref{thm:stability_obs} can be interpreted as describing the stability of $\postmeas$ in the presence of experimental sources of error that impact the ability to determine $\obsmeas$ to a specified level of accuracy. 
On the other hand, even when $\priormeas$ and $\obsmeas$ are given and considered to be ``error free'',  we must still determine $\pfpriormeas$ in order to determine or sample from $\postmeas$.
This is often done numerically, e.g., as shown in the following section, which introduces errors into the conditional densities.
Specifically, if $\widehat{\pfpriordens}$ denotes the computed push-forward density, then the conditional densities are defined, in the sense of the disintegration theorem, to be
$$
	\frac{\widehat{dP_q^{\text{prior}}}}{d\mu_{\pspace,q}(\lambda)} = \frac{\priordens(\lambda)}{\widehat{\pfpriordens}(Q(\lambda))}.
$$
In Section~\ref{sec:comp} we give a precise bound on the total variation error in the posterior in terms of the total variation error in the push-forward of the prior.
\section{Computational Considerations}\label{sec:comp}
We note that both $\priordens$ and $\obsdens$ are assumed to be given, so we only need to approximate the push-forward probability 
density, $\pfpriordens$, induced by the prior and the model. 
In general, we cannot expect to derive an analytical form for the push-forward of the prior, so we require 
a means to numerically approximate this probability density function.  As previously mentioned, our focus in this paper 
is on the exposition of this new approach, so we utilize straightforward Monte Carlo sampling to propagate $\priordens$
and $\postdens$ through the model and a kernel density approach~\cite{wand1994multivariate} to approximate $\pfpriordens$ and $\pfpostdens$.
This procedure is summarized in Algorithm~\ref{alg:PF}.
\begin{algorithm}
\caption{Computing the Push-Forward Density Induced by the Prior and Model}\label{alg:PF}
 Given a set of samples from the prior density: $\lambda_i, \quad i=1,\ldots,M$\;
 Evaluate the model and compute the QoIs: $q_i = Q(\lambda_i)$\; 
Use the set of QoIs and use a standard technique, such as a kernel density approximation to estimate $\pfpriordens(q)$\;
\end{algorithm}
We emphasize that Alg.~\ref{alg:PF} is simply a forward propagation of uncertainty using Monte Carlo sampling and an approach for density estimation.
While Monte Carlo sampling is easy to implement and can effectively utilize parallel model evaluations, the rate of convergence is relatively low (${\cal O}(N^{-1/2})$).
More advanced methods based on global polynomial approximations computed using stochastic spectral methods \cite{GhanemSpanos, XiuKarniadakis, WanKarniadakis, MaitreGhanem1, GhanemRedhorse}, sparse grid collocation methods \cite{gerstner03,hegland2003,bungartz04,Ma:2009:AHS:1514432.1514547}, and Gaussian process models~\cite{rasmussen2006}
converge much faster (under certain assumptions) and can also be used estimate the push-forward of the prior.
The utilization of these methods within the consistent Bayesian framework will be thoroughly investigated in future work.

Clearly, the accuracy of the posterior density depends on the accuracy of the approximation
of the push-forward of the prior.  In this paper, we utilize the standard Gaussian kernel density estimation schemes which are known to converge ${\cal O}(N^{-4/(4+d)})$ in the mean-squared error~\cite{Terrel_S_JSTOR_1992} and ${\cal O}(N^{-2/(4+d)})$ in the L$_1$-error~\cite{Devroye85}.
We use $\widehat{\pfpriordens}$ to denote the approximation of the push-forward of the prior density and we use $\widehat{\postdens}$ to denote the corresponding approximation of $\postdens$, i.e.,
\[\widehat{\postdens}(\lambda) = \priordens(\lambda) \frac{\obsdens(Q(\lambda))}{\widehat{\pfpriordens}(Q(\lambda))}.\]
We make the following assumption regarding the compatibility of this approximation of the push-forward of the prior and the observed measure.
\begin{assumption}\label{assump:approxdom}
There exists a constant $C>0$ such that $\obsdens(q)\leq C\widehat{\pfpriordens(q)}$ for a.e.~$q\in\dspace$.
\end{assumption}

We can now prove the following result concerning the accuracy of the posterior in terms of the accuracy of the approximation of the push-forward of the prior.
\begin{theorem}\label{thm:acc_posterior}
If the approximation of the push-forward of the prior satisfies Assumption~\ref{assump:approxdom}, then the error in the posterior can be bounded by,
\[ d_{TV}(\postmeas,\widehat{\postmeas}) \leq C d_{TV}(\pfpriormeas,\widehat{\pfpriormeas}).\]
\end{theorem}
\begin{proof}
Assume $\widehat{\pfpriormeas}$ satisfies Assumption~\ref{assump:approxdom}.  Then we have,
\begin{eqnarray*}
d_{TV}(\postmeas,\widehat{\postmeas}) &=& \int_{\dspace} \bigg(  \int_{\pspace\cap Q^{-1}(q)} \, \abs{\frac{\priordens(\lambda)}{\pfpriordens(q)} - \frac{\priordens(\lambda)}{\widehat{\pfpriordens(q)}}} \obsdens(Q(\lambda)) d\mu_{\pspace,q}(\lambda) \bigg) \, d\dmeas(q) \\
&=& \int_{\dspace} \bigg(  \int_{\pspace\cap Q^{-1}(q)} \, \priordens(\lambda) d\mu_{\pspace,q}(\lambda) \bigg) \abs{\frac{1}{\pfpriordens(q)} - \frac{1}{\widehat{\pfpriordens(q)}}} \obsdens(q)\, d\dmeas(q) \\
&=& \int_{\dspace} \, \pfpriordens(q) \, \abs{\frac{1}{\pfpriordens(q)} - \frac{1}{\widehat{\pfpriordens(q)}}} \obsdens(q)\, d\dmeas(q) \\
&=& \int_{\dspace} \, \frac{\obsdens(q)}{\widehat{\pfpriordens(q)}} \, \abs{\widehat{\pfpriordens(q)} - \pfpriordens(q)} \, d\dmeas(q) \\
&\leq & C d_{TV}(\pfpriormeas,\widehat{\pfpriormeas}),
\end{eqnarray*}
where the last step follows from Assumption~\ref{assump:approxdom}.
\end{proof}
%

In our approach, we first perform a forward propagation of the prior density to obtain an approximation of $\pfpriordens$.  Given this push-forward
density, we can directly interrogate the posterior density for any $\lambda \in \pspace$ given $Q(\lambda)$.
To simplify the discussion, we rewrite our posterior density~\eqref{eq:postpdf} as
\[\postdens(\lambda) = \priordens(\lambda)\ r(Q(\lambda)), \quad \text{ where } \quad r(Q(\lambda)) = \frac{\obsdens(Q(\lambda))}{\pfpriordens(Q(\lambda))}.\]
We use $r(Q(\lambda))$ instead of simply $r(\lambda)$ to remind the reader that for a given $\lambda \in \pspace$, we need to compute $Q(\lambda)$ which requires a model evaluation in order to evaluate the posterior.
Throughout this paper, we leverage the fact that we have already evaluated the model for the samples generated from the prior that we used in Alg.~\ref{alg:PF}.
In fact, many standard calculations involving the posterior only involve integrals of $r(Q(\lambda))$ with respect to the prior.
For example, the integral of the posterior is given by
\[\text{I}(\postdens) = \int_{\pspace} \postdens(\lambda) \ d\pmeas = \int_{\pspace} \priordens(\lambda) r(Q(\lambda)) \ d\pmeas = \int_\pspace r(Q(\lambda)) \ d\priormeas,\]
and the Kullback-Liebler (KL) divergence is given by,
\[\text{KL}(\priordens : \postdens) = \int_{\pspace} \postdens(\lambda) \log{\left(\frac{\postdens(\lambda)}{\priordens(\lambda)}\right)} \ d\pmeas = \int_\pspace r(Q(\lambda))\log{r(Q(\lambda))}\ d\priormeas.\]

A common goal in stochastic inversion is to generate a set of samples from the posterior which can then
be used to characterize the posterior density or to approximate the push-forward of the posterior
and make predictions on QoI that cannot be measured or observed.
We can easily apply rejection sampling to the samples from prior (generated in Alg.~\ref{alg:PF}) using only
$r(Q(\lambda))$.
This algorithm is summarized in Alg.~\ref{alg:GS}.



\begin{algorithm}
\caption{Generating samples from the posterior using rejection}\label{alg:GS}
   Given the $\pfpriordens(q)$ from Algorithm~\ref{alg:PF}\;
   Select a set of samples in $\pspace$: $\lambda_p, \quad p=1,\ldots,P$\;
   Evaluate the model and compute the QoIs: $q_p = Q(\lambda_p)$\; 
   Compute $r(Q(\lambda))$ at each $\lambda_p$\;
   Estimate $M = \max_{\pspace} r(Q(\lambda))$\;
   \For{$p=1,\ldots,P$}{
     Generate a random number, $\xi_p$, from a uniform distribution on $[0,1]$\;
     Compute the ratio: $\eta_p = r(Q(\lambda))/M$\;
     \eIf{$\eta_p > \xi_p$} {
       Accept $\lambda_p$\;
     }{
       Reject $\lambda_p$\;
     }
   } 
\end{algorithm}
We use $r(Q(\lambda))$ in Alg.~\ref{alg:GS} to emphasize the fact that applying rejection to generate samples from the posterior using
samples generated from the prior in $\pspace$ is equivalent to applying rejection to generate samples from the observed density using samples generated from the push-forward of the prior in $\dspace$.
Since we assume $\dspace$ is relatively low-dimensional, i.e., we only have a small number of QoI, we do not observe significant degradation of the rejection sampling if the dimension of $\pspace$ is large.
Obviously if $\pfpriordens$ and $\obsdens$ are significantly different, i.e., if the data is very informative,
then the percentage of accepted samples may be quite low.
Additional samples can be generated using a variety of techniques, e.g., importance sampling or MCMC,
and we will explore efficient utilization of these techniques for the consistent Bayesian approach
in future work.




\section{Numerical examples}\label{sec:numerics}

In this section we present some numerical results to verify the theoretical results presented in
Section~\ref{sec:Bayes}.
We start with a simple parameterized nonlinear system with two parameters and then consider
a two-dimensional discontinuous function.
We then consider a 100-dimensional example based on a finite element
discretization of a model for single-phase incompressible flow in porous media
with a Karhunen-Loeve expansion for the log-transformed permeability field.
We conclude this section with an analytical example to demonstrate the convergence of the
posterior as the approximation of the push-forward of the prior improves, e.g., as more samples
are evaluated.
We investigate this rate of convergence in terms of the dimension of the parameter space and in terms of the
dimension of the observation space.

\subsection{A Parameterized Nonlinear System}\label{subsec:ex2}

Consider the following parameterized nonlinear system of equations was introduced in \cite{BE1}:
\begin{eqnarray*}
\lambda_1 x_1^2+ x_2^2 &=& 1 \\
x_1^2- \lambda_2 x_2^2 &=& 1
\end{eqnarray*}

The first QoI is the second component, i.e., $q(\lambda) = x_2(\lambda)$.
The parameter ranges are given by $\lambda_1 \in [0.79,0.99]$ and $\lambda_2 \in [1-4.5\sqrt{0.1},1+4.5\sqrt{0.1}]$ 
which are chosen as in \cite{BE1} to induce 
an interesting variation in the QoI.
The QoI as a function of the parameters over the given
ranges can be found in \cite{BE1}.
We assume that the 
observed distribution on the QoI is a normal distribution 
with mean 0.3 a standard deviation of 0.025.
We consider two different prior probability densities on $\pspace$, uniform and Beta$(2,5)$, and show that despite the 
fact that they lead to different posteriors, the push-forward densities induced by the posteriors both match the observed density.

We follow Algorithm~\ref{alg:PF} and generate 10,000 samples from each prior and use
a kernel density estimator to approximate the resulting push-forward of the prior densities.
Using these approximations of $\pfpostdens(q)$,
we follow Algorithm~\ref{alg:GS} to
generate a set of samples from each posterior density.
As described in Section~\ref{sec:comp}, we can use the values of $r(Q(\lambda))$ evaluated at the samples generated from the prior to approximate the integral of the posterior as well as the KL-divergence from the prior to the posterior:
\[ \text{I}(\postdens) \approx 0.9993, \quad \text{KL}(\priordens : \postdens) \approx 1.1344.\]

\begin{figure}[ht]
 \begin{center}
\includegraphics[width=0.32\textwidth]{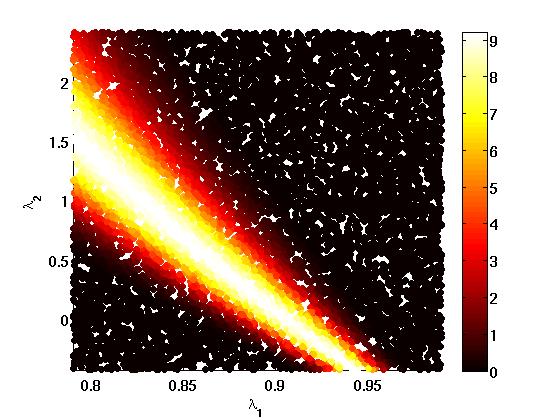}
\includegraphics[width=0.32\textwidth]{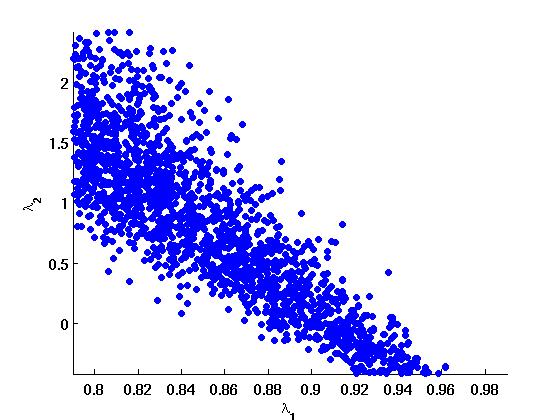}
\includegraphics[width=0.32\textwidth]{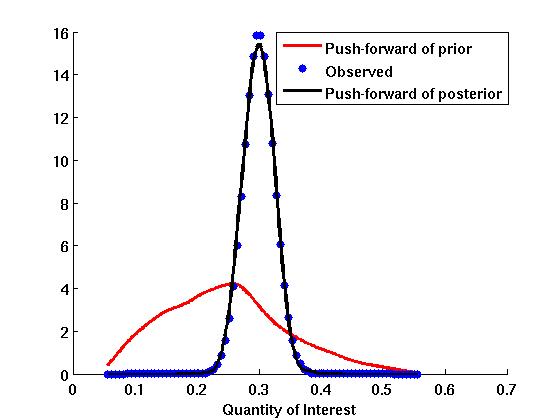}
 \end{center}
 \caption{Approximation of the posterior density obtained using consistent Bayesian inference with a uniform prior (left), 
 set of samples from the posterior (middle) and a comparison of the observed density with the push forward densities obtained by propagating the prior and the posterior densities through the forward model}  
 \label{fig:ex2_unif}
 \end{figure}

In Figure~\ref{fig:ex2_unif}, we plot the posterior density evaluated at the samples from the uniform prior the set of the
samples that we generate from each posterior using the accept/reject criteria of Algorithm~\ref{alg:GS}, and the comparison of the push-forward of the prior and the posterior with the observed density.

\begin{figure}[ht]
 \begin{center}
\includegraphics[width=0.32\textwidth]{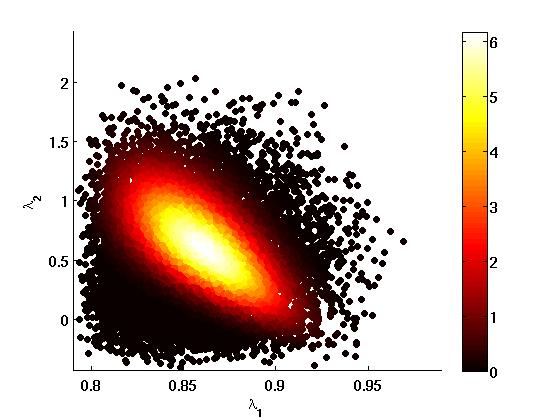}
\includegraphics[width=0.32\textwidth]{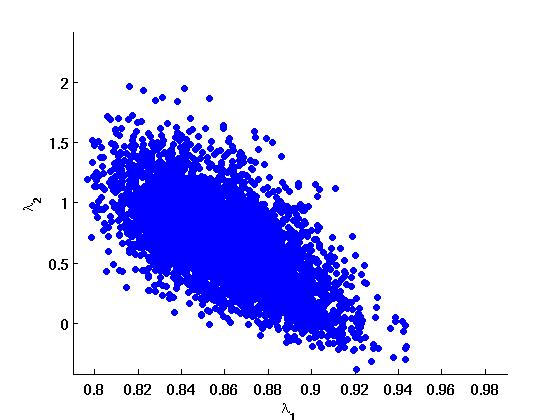}
\includegraphics[width=0.32\textwidth]{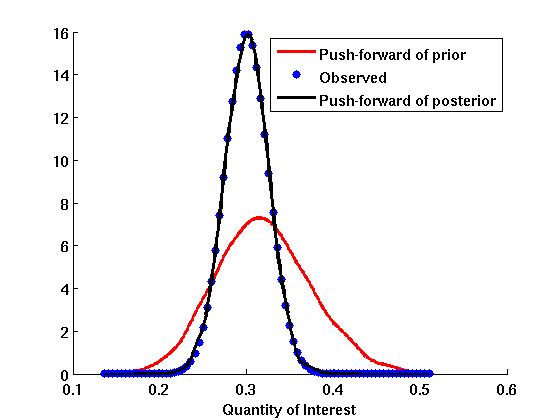}
 \end{center}
 \caption{Approximation of the posterior density obtained using consistent Bayesian inference with a Beta(2,5) prior (left), 
 set of samples from the posterior (middle) and a comparison of the observed density with the push forward densities obtained by propagating the prior and the posterior densities through the forward model}  
 \label{fig:ex2_beta}
 \end{figure}

 In Figure~\ref{fig:ex2_beta}, we plot the approximations of the posterior density evaluated at the samples from the Beta$(2,5)$ prior, the set of the
samples that we generate from each posterior using the accept/reject criteria of Algorithm~\ref{alg:GS}, and the comparison of the push-forward of the prior and the posterior with the observed density.
As before, we use the values of $r(Q(\lambda))$ evaluated at the samples generated from the prior to approximate the integral of the posterior as well as the KL-divergence from the prior to the posterior:
\[ \text{I}(\postdens) \approx 1.0106, \quad \text{KL}(\priordens : \postdens) \approx 0.4399.\]

Clearly, the push-forward of the two priors differ significantly from the observed density and from each other, but the push-forward of the two posteriors
both agree quite well with the observed density.

\subsection{A Piecewise Smooth Example}\label{subsec:ex3}
To demonstrate the ability of our approach to generate a consistent posterior for piecewise smooth functions, we consider the d-dimensional function introduced in \cite{Jakeman2013790}:
\[ q(x) = 
\begin{cases}
q_1(x)-2, & 3x_1+2x_2 \geq 0 \text{ and } -x_1+0.3x_2 < 0, \\
2q_2(x), & 3x_1+2x_2 \geq 0 \text{ and } -x_1+0.3x_2 \geq 0, \\
2q_1(x)+4, & (x_1+1)^2 + (x_2+1)^2 < 0.95^2 \text{ and } d=2,\\
q_1(x), & \text{otherwise}
\end{cases}, \]
where
\[ q_1(x) = \exp\left( -\sum_{i=1}^d x_i^2\right) - x_1^3 - x_2^3, \quad q_2(x) = 1+q_1(x)+\frac{1}{4d}\sum_{i=1}^d x_i^2.\]
For visualization purposes, we focus on $d=2$.  We direct the interested reader to \cite{Jakeman2013790} for a visualization of this response surface.
This particular example demonstrates one of the challenges in performing inversion with discontinuous functions.  
Namely, with $\pspace = [-1,1]^2$, $Q(\pspace)$ is a disconnected (but compact) region in $\mathbb{R}$.  
As previously mentioned, our approach assumes that the observed probability density is defined over 
$\dspace = Q(\pspace)$ which may not be known {\em a priori}.
Fortunately, it is fairly easy to determine when this assumption is violated (the integral of the
posterior will not be 1) and since the main computational expense in the consistent Bayesian approach is
 computing the push-forward of the prior, we can explore utilizing different observed densities with only 
a small computational effort.  
We assume a uniform prior on $\pspace$ and we define $\obsdens \sim N(-2.0,0.25^2)$.
We use the values of $r(Q(\lambda))$ evaluated at the samples generated from the prior to approximate the integral of the posterior as well as the KL-divergence from the prior to the posterior:
\[ \text{I}(\postdens) \approx 0.9998, \quad \text{KL}(\priordens : \postdens) \approx 1.6014.\]

\begin{figure}[ht]
 \begin{center}
   \includegraphics[width=0.32\textwidth]{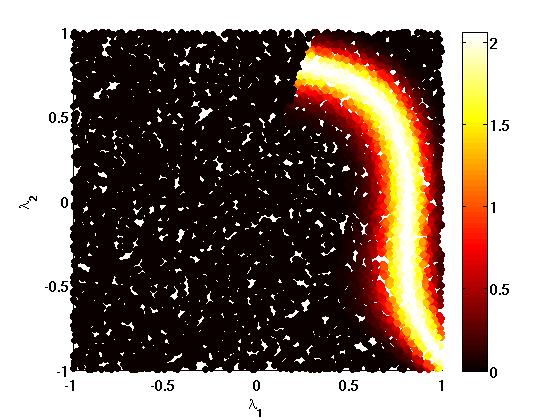}
   \includegraphics[width=0.32\textwidth]{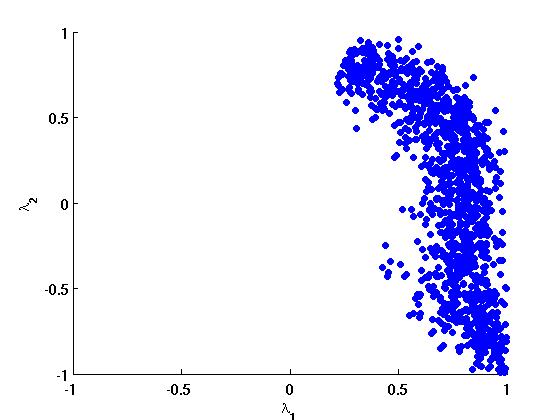}
   \includegraphics[width=0.32\textwidth]{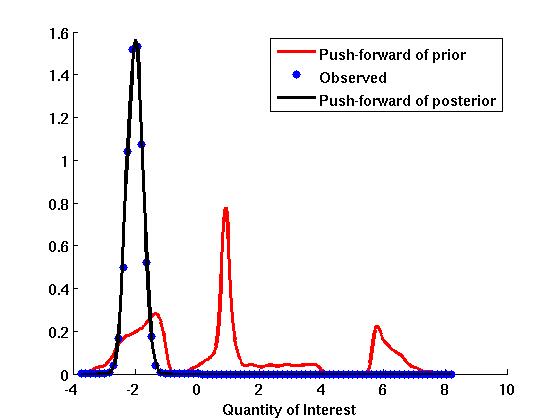}
 \end{center}
 \caption{ Approximation of the posterior density obtained using consistent Bayesian inference on the piecewise smooth function with a uniform prior (left), 
 set of samples from the posterior (middle) and a comparison of the observed density with the push forward densities obtained by propagating the prior and the posterior densities through the forward model.}
 \label{fig:disc_images}
 \end{figure}

As in the previous section, following Algorithm~\ref{alg:PF}, we use 10,000 samples (generated from the prior) to compute the push-forward of the prior.
In Figure~\ref{fig:disc_images}, we plot the posterior probability density evaluated at the samples from the prior, the subset of the samples used to compute the push-forward of the prior that passed the accept/reject
criteria given in Algorithm~\ref{alg:GS}, and the probability densities, 
$\obsdens$, $\pfpriordens$ and $\pfpostdens$, on $\dspace$.   
We see that the push-forward of the prior indicates that $\dspace$ is a disconnected region in $\mathbb{R}$.
The observed density satisfies the assumptions described earlier, so the consistent Bayesian approach
gives a posterior probability density on $\pspace$, which, when propagated through the model, agrees quite
well with the observed density.

\subsection{A Higher-dimensional PDE-based Example}\label{subsec:ex4}
Consider a single-phase incompressible flow model:
\begin{equation}\label{eq:porous}
\begin{cases}
-\nabla \cdot (K(\lambda) \nabla p) = 0, & x\in\Omega = (0,1)^2,\\
p = 1, & x=0, \\
p = 0, & x=1, \\
K\nabla p \cdot \mathbf{n} = 0, & y=0 \text{ and } y=1.
\end{cases}
\end{equation}
Here, $p$ is the pressure field and $K$ is the permeability field which we assume is a scalar field given by a Karhunen-Lo\'eve expansion of the log transformation, $Y = \log{K}$, with
\[Y(\lambda) = \overline{Y} + \sum_{i=1}^\infty \xi_i(\lambda)\sqrt{\eta_i}f_i(x,y),\]
where $\overline{Y}$ is the mean field and $\xi_i$ are mutually uncorrelated random variables with zero mean and unit variance \cite{ganis2008stochastic,wheeler2011multiscale}.
The eigenvalues, $\eta_i$, and eigenfunctions, $f_i$, are computed using an assumed functional form for the covariance matrix \cite{zhang2004efficient,Schwab2006100}.
We assume a correlation length of $0.01$ in each spatial direction and truncate the expansion at 100 terms.
This choice of truncation is purely for the sake of demonstration.
In practice, the expansion is truncated once a sufficient fraction of the energy in the eigenvalues is retained~\cite{zhang2004efficient,ganis2008stochastic}.
To approximate solutions to the PDE in Eq.~\eqref{eq:porous} we use a finite element discretization with continuous piecewise bilinear basis functions defined on a uniform ($50\times 50$) spatial grid.

Our quantity of interest is the pressure at $(0.0540,0.5487)$.
The prior is a multivariate standard normal density $\priordens \sim N({\mathbf 0},{\mathbf I})$ where ${\mathbf I}$ is the standard identity matrix.
Since the parameter space is 100-dimensional and we use a standard multi-variate normal for the prior, numerical integration of either the prior or the posterior with respect to the volume measure is quite challenging.
However, we only need to integrate functions of $r(Q(\lambda))$ with respect to the prior.
We generate 10,000 samples from the prior and evaluate the PDE model for each of these realizations.
We use a standard KDE to approximate the push-forward of the prior in the 1-dimensional output space.
We assume the observed density on the QoI is given by $\obsdens \sim N(0.7,1.0\text{E-4})$.
We use the values of $r(Q(\lambda))$ evaluated at the samples generated from the prior to approximate the integral of the posterior as well as the KL-divergence from the prior to the posterior:
\[ \text{I}(\postdens) \approx 0.9968, \quad \text{KL}(\priordens : \postdens) \approx 0.7372.\]

We apply rejection sampling to select a subset of the samples from the prior for the posterior.
This gives 2781 samples from the posterior ($\approx 28\%$ acceptance).
We again emphasize that we are only using rejection sampling to select samples from the observed density using the samples generated from the push-forward of the prior in $\dspace$.
Thus, even for this moderately high-dimensional parameter space, we see a reasonable acceptance rate.

\begin{figure}[ht]
\begin{center}
\includegraphics[width=0.4\textwidth]{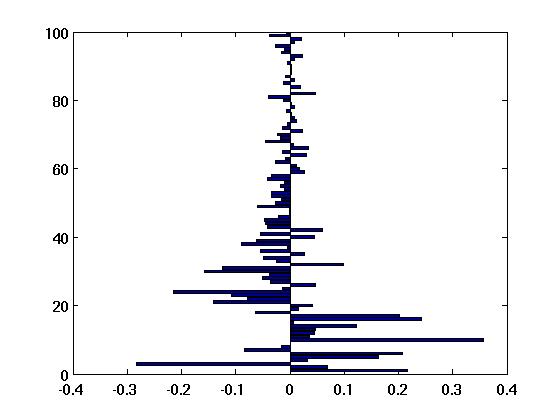}
\includegraphics[width=0.4\textwidth]{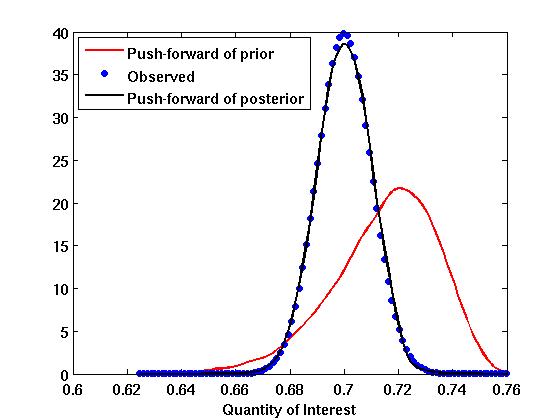}
\end{center}
\caption{ The means of the accepted samples from the posterior (left) and a comparison of the observed density with the push-forward of the prior and posterior densities (right).}
\label{fig:100D_images}
\end{figure}

Describing a non-parametric density in a 100-dimensional space is challenging.
On the left side of Figure~\ref{fig:100D_images}, we plot the mean of the accepted samples in each dimension.
Recall that the mean of the prior was zero is every direction.
These shifts in the means clearly indicate that the data has informed the model is some way, but precisely
characterizing the posterior is difficult since it is no longer Gaussian.
However, we know that the push-forward of the posterior should match the observations and this is shown on the right side of Figure~\ref{fig:100D_images}.

\subsection{Numerical convergence}\label{sec:numerical-convergence}
Here, we investigate the effect of increasing parameter dimension on the accuracy of numerically constructed posterior densities.
Consider the following model
\begin{align}\label{eq:chi2-quadratic}
Q(\lambda) = (\lambda-\mu)^TC^{-1}(\lambda-\mu).
\end{align}
If the prior density $\priordens$ is $d$-dimensional Gaussian with mean $mu$ and covariance $C$, then the push-forward probability density of~\eqref{eq:chi2-quadratic} is a univariate Chi-squared distribution with parameter $d$, i.e.,
\begin{align}\label{eq:chi-squared-sensity}
\pfpriordens(Q(\lambda))=\pi_{\chi^2}(Q(\lambda))={\frac {1}{2^{\frac {d}{2}}\Gamma \left({\frac {d}{2}}\right)}}\;Q(\lambda)^{{\frac {d}{2}}-1}e^{-{\frac {Q(\lambda)}{2}}}.
\end{align}
In the following, we set $\mu=(0,\ldots,0)^T$ and define a random matrix $A\in\mathbb{R}^{d\times d}$ with entries drawn randomly from the standard normal distribution and set $C=A^TA$.

\begin{figure}[ht]
\begin{center}
\includegraphics[width=0.45\textwidth]{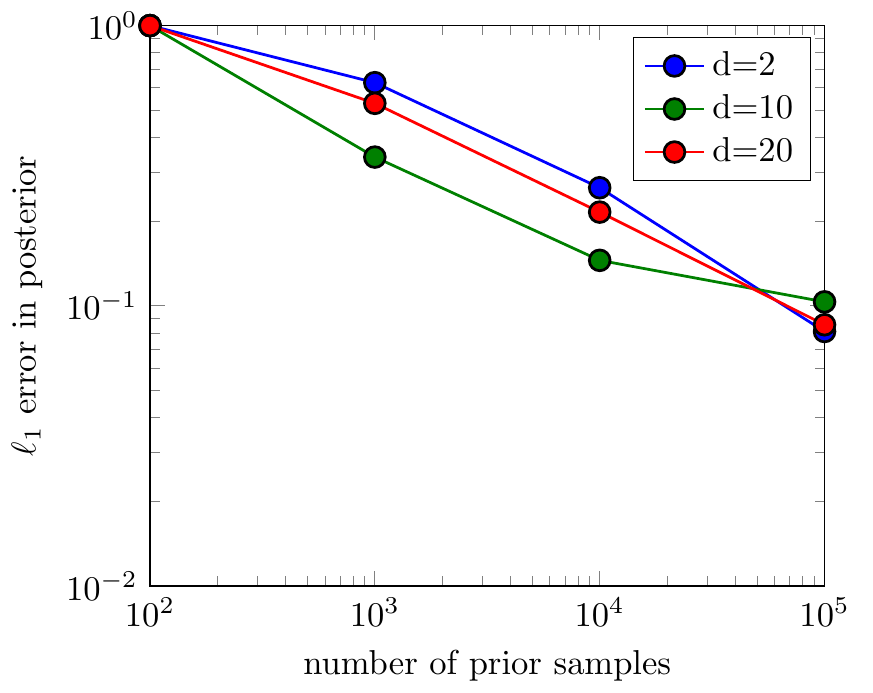}
\includegraphics[width=0.45\textwidth]{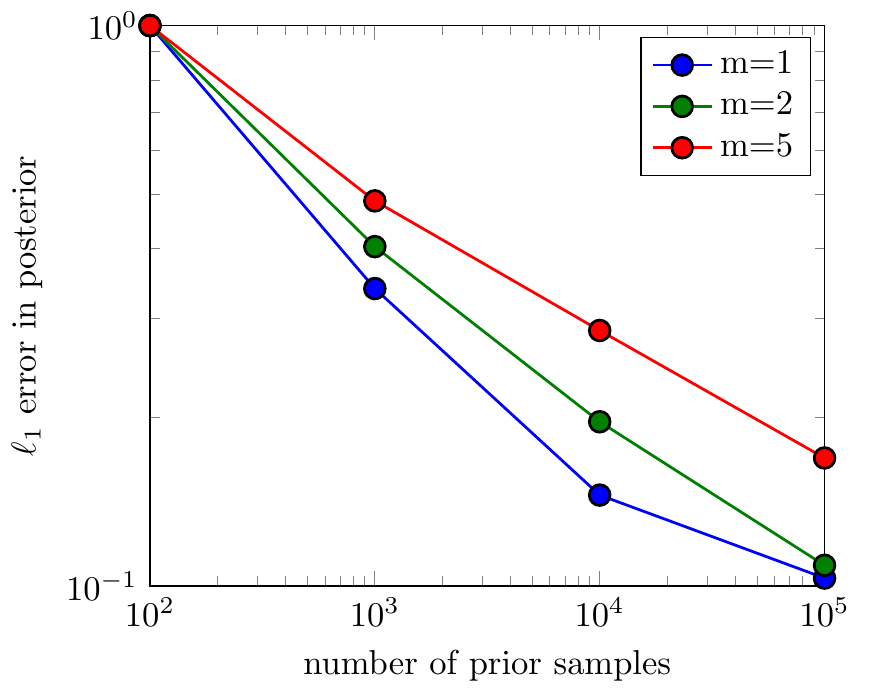}
\end{center}
\caption{Convergence of in the error of the posterior when the $\widehat{\pfpriordens}$ is computed using a GKDE. Convergence is shown for increasing number of model evaluations, for increasing dimension (left) and numbers of QoI (right). To facilitate easier comparison each curve represents the decrease in error relative to the error when 100 samples are used to build the GKDE.}
\label{fig:chi2-error-convergence}
\end{figure}

Using the exact expression for the push-forward of the prior, we compute the posterior density analytically for a given observed density. 
In the following, we assume the observed density to be uniformly distributed over the range $[a,b]$, where we vary the parameters $a$ and $b$ as we change the dimensionality $d$ of the prior. 
To ensure a fair comparison, we choose $a$ and $b$ so that the measure of the observed density with respect to the push-forward measure remains the same for each dimension considered. 
Specifically, we choose $a$ such that
$$
	\int_{-\infty,a}\pi_{\chi^2}(Q(\lambda))\, d\dmeas(q)=\frac{2}{5},
$$
and $b$ is chosen such that
$$
	\int_{-\infty,b}\pi_{\chi^2}(Q(\lambda))\, d\dmeas(q)=\frac{3}{5}.
$$
In Figure~\ref{fig:chi2-error-convergence}, we plot the convergence of the error of the posterior using the approximate push forward $\widehat{\pfpriordens}$, computed using Gaussian kernel density estimation (GKDE), against the number of samples $N$ used to build the GKDE.
Here we set the bandwidth of the KDE based upon the optimal asymptotic width and as theory suggests~\cite{Terrel_S_JSTOR_1992} the convergence rate is $O(N^{-\frac{2}{5}})$, which is independent of the {\em parameter} dimension.

The curves shown in Figure~\ref{fig:chi2-error-convergence} reflect the median error taken over 20 random experiments. 
We use a discrete approximation of the $L_1$-norm, evaluated using $N=10,000$ samples from the prior, to measure the error in the GKDE (i.e., we approximate the total variation metric). 
Specifically,
$$\lVert\widehat{\postdens}(\lambda)-\postdens(\lambda)\rVert_1\approx \frac{1}{N}\sum_{i=1}^{N} \left|\frac{\obsdens(Q(\lambda_i))}{\pfpriordens(Q(\lambda_i))}-\frac{\obsdens(Q(\lambda_i))}{\widehat{\pfpriordens}(Q(\lambda_i))}\right| $$
Again, we have leveraged the fact that we can compute integrals with respect to the prior measure using the samples used to compute the push-forward of the prior.

Now, we investigate the effect of increasing dimension of the output space $\dspace$ on the accuracy of numerically constructed posteriors.
Considered a modified version of the previous model that has $m$ QoI. 
In contrast to the previous example, we set the prior to be a multivariate Gaussian with a block diagonal covariance with $m$ blocks. 
Each block $C_i\in\mathbb{R}^{d/m\times d/m}$, $i=1,\ldots,m$ is generated from random matrices $A_i$, and we set the mean $\mu_i$ associated with each block to be zero. 
With this new definition of the prior, we now set the $i$-th QoI of the model to be
\begin{align}\label{eq:chi2-quadratic-block}
Q_i(\lambda) = (\lambda-\mu_i)^TC_i^{-1}(\lambda-\mu_i).
\end{align}
The exact push-forward of the prior for each QoI is again the Chi-squared distribution, and the joint density  is simply the tensor product of the one-dimensional marginals. 
We remark that although each QoI is only dependent on $m/d$ variables, when inferring the parameters using the joint density of the push-forward, each QoI will contribute to the form of the posterior. 
Again, we use a uniform observed density, this time over $[a,b]^m$, where we now choose $a$ such that
$$
	\int_{-\infty,a}\pi_{\chi^2}(Q(\lambda))\, d\dmeas(q)=\frac{1}{2}-(\frac{1}{5})^{1/m},
$$
and $b$ is chosen such that
$$
	 \int_{-\infty,b}\pi_{\chi^2}(Q(\lambda))\, d\dmeas(q)=\frac{1}{2}+(\frac{1}{5})^{1/m}.
$$
These choices of $a$ and $b$ ensure that the push-forward measure remains the same for each dimension $m$ of the output space considered.

In Figure~\ref{fig:chi2-error-convergence}, we plot the convergence of the error of GKDE based posteriors for increasing numbers, $m$, of the QoI. 
Although the convergence rate of the error in the posterior is independent of parameter dimension $d$, it is heavily dependent on the dimensionality $m$ of the ouput space. Specifically the convergence rate is $O(N^{-\frac{2}{m+4}})$, using the asymptotically optimal bandwidth~\cite{Devroye85,Terrel_S_JSTOR_1992}. There are many alternative density estimation techniques that can be used to construct $\pfpriordens$ and in some cases these alternative methods can be much more effective than the form of the GKDE we used here~\cite{Royset_W_EJOR_2015,Terrel_S_JSTOR_1992}. The exploration of alternative density estimation techniques will be a topic of future work.


\section{Comparison with Existing Approaches}\label{sec:compare}
As mentioned in Section~\ref{sec:intro}, our approach was inspired by both the statistical Bayesian approach, which is perhaps
the most commonly used approach to perform stochastic inference, and the recently developed measure-theoretic approach for
solving stochastic inference problems \cite{BE1,BE2,BE3}.  In this section, we briefly compare our approach with these existing
approaches.

\subsection{The Statistical Bayesian Approach}\label{subsec:compstat}

In this section, we present the traditional Bayesian approach to the inverse problem of finding the model input $\lambda$ given observations (or functionals) of the solution of the model. Here we will refer to this approach as Statistical Bayesian inference, reflecting the nature of the numerical sampling methods, such as Markov Chain Monte Carlo (MCMC), typically used to solve the inverse problem.

Following \cite{Stuart_AN_2010} we assume that observational data, $q\in \dspace$ are noisy outputs to a governing forward problem, that is
\begin{equation}
q=Q(\lambda)+\eta,\label{eq:stat-bayes-additive-obs-error}
\end{equation}
where $\eta$ is a random variable with mean zero. Assuming that $\eta\in\mathbb{R}^m$ is a random variable with probability density $\rho$, the probability of the data given the random variables is given by the likelihood function
$$\pi(q|\lambda)=\rho(q-Q(\lambda)).$$
The likelihood explicitly represents the probability that a given set of parameters might give rise to the observed data.

Following this Bayesian formulation, the posterior density is
\begin{equation}\label{eq:statbayes}
\postdenssbayes(\lambda|q) = \priordens(\lambda) \frac{\pi(q|\lambda)}{\int_{\Lambda} \pi(q|\lambda) \priordens d\pmeas}.
\end{equation}
However, since the denominator is a normalizing constant not necessary for implementing sampling methods such as MCMC, this is often re-written, for convenience, as
\[\postdenssbayes(\lambda|q) \propto \priordens(\lambda) \pi(q|\lambda).\]
As noted by \cite{Stuart_AN_2010}, the above relationship can be used to
express the relationship of the posterior measure $\tilde{P}^{\mathrm{prior}}_\pspace$ and the prior measure $\priormeas$ as the Radon Nikodym derivative
$$
\frac{d\tilde{P}^{\mathrm{post}}_\pspace}{d\priormeas}\propto \rho(q-Q(\lambda))
$$

The probability density of $\eta$, and thus likelihood function, is a modeling assumption. It is quite common, but certainly not necessary, to assume that the difference between the 
model prediction and the  data is given by a Gaussian random variable, that is $\eta \sim N(0,\Sigma)$
which yields
$$\postdenssbayes(\lambda|q) \propto \priordens(\lambda) \exp(-\frac{1}{2}q^T\Sigma^{-1} q)$$
where the covariance matrix, $\Sigma$, is given based on the noise in the measurements.

The assumption of Gaussian noise is by no means the only choice, or indeed always the best one. The model $Q(\lambda)$ is often only an approximation of the process it is modelling. In this case it may be more appropriate to replace ~\eqref{eq:stat-bayes-additive-obs-error} with $q=Q(\lambda)+\delta+\eta$, where $\delta\in\mathbb{R}^m$ represents the  discrepancy between the model and the true noiseless processing being simulated. Alternative discrepancy models can also be used as an alternative to the aforementioned additive discrepancy\footnote{Recently there has been some success embedding physics-based  discrepancy within a model~\cite{Sargsyan_NG_IJCK_2015}.}.

The transformation of the stochastic inverse problem to an inference problem makes a direct comparison challenging.  In Section~\ref{subsec:ex1}, we compare the posterior distributions using the consistent and statistical Bayesian approaches for a simple model problem and show that the two are identical only for special choices of the prior and the model, and here we enumerate some of the key conceptual differences of the two approaches.

Statistical and consistent Bayesian inference both quantify the
relative probability that a prediction made by a mathematical model will lie in certain specified regimes. However, the assumptions made by these two approaches when addressing this question differ, so it is perhaps no surprise that the push-forward probabilities differ as well. 
Statistical Bayesian inference assumes that the difference between the model and the observational data is a random variable with density $\rho$. 
However, when the posterior density obtained is pushed forward through the model, the resulting density is not typically $\rho$. 
In contrast, consistent Bayesian inference does not explicitly assume an error model, but rather assumes (or is given) an observed density on the QoI and finds a posterior such that the push-forward of the posterior through the model 
exactly reproduces the density on the observations.

In addition, consistent Bayesian inference requires global knowledge of the forward map over the non-zero support of the prior while statistical Bayesian inference only requires local knowledge over the intersection of the non-zero support of the prior and likelihood function.
The advantage of having global knowledge of the map is that once the forward map is known, and new observational data becomes available, the posterior distribution can be updated without additional evaluations of the forward model. 
The disadvantage of requiring global knowledge of the forward map is that typically this knowledge is gained through exhaustive sampling throughout the parameter space which can be computationally expensive.

The statistical Bayesian and deterministic optimization perspectives are linked through the posterior probability density and the maximum a posteriori (MAP)
estimator, which maximizes the posterior density (see e.g., \cite{Stuart_AN_2010,Bui-ThanhGhattasMartinEtAl13}).
The connections between the consistent Bayesian approach and deterministic optimization have not thoroughly explored and are the subject of current research.


In general, these two approaches introduce stochasticity in different ways, solve different problems, give different posteriors and make different predictions.
The choice on which to use should be based on the appropriate problem formulation which will be problem dependent.
For example, if the goal is to find a posterior that produces a specific push-forward (posterior predictive), then the consistent Bayesian approach would be preferable.
In addition, if the goal in solving the stochastic inverse problem is to subsequently use the posterior to predict QoI that cannot be observed directly, which corresponds to a puch-forward of the posterior to the new observation space, then it may be preferable to use a posterior with known posterior predictive for the QoI that can be measured.
On the other hand, there are certainly scenarios where the ``classical'' Bayesian approach may be preferable.
For example, if there process used to produce the prior were more reliable than the process used to measure the data, then the consistent Bayesian approach may not be the proper choice since it completely ignores the prior in the directions informed by the data.  The ``classical'' Bayesian approach always include some influence from the prior, even in the directions informed by the data.


\subsection{An Illustrative Example}\label{subsec:ex1}

In this section, we describe a simple example to compare the statistical and consistent Bayesian approaches.
Consider the simple mapping from a 1-dimensional parameter space, $\pspace = [-1,1]$,
to a single QoI.  The mapping is given by
\[ q(\lambda) = \lambda^p, \quad p=1,3,5\ldots.\]
The parameter, $p$, is not uncertain and will be used to vary the nonlinearity of the map.
We assume a uniform prior and that the uncertainty on the observations is
a truncated normal density, $\obsdens(q) \sim N(\hat{q},\sigma^2)$,
where $\hat{q}$ and $\sigma$ are also assumed to be given.
While this density is truncated by $q(\pspace)$ (and is also re-normalized), we use the above form of the density for simplicity. 
Our assumptions on the observed density for the QoI are consistent with an 
observed value of $\hat{q}$ and assuming $\eta \sim N(0,\sigma^2)$ in the statistical Bayesian approach. 
Since the prior is uniform, the two posterior measures/densities are {\em identical} if $p=1$, i.e., if the map is linear.
This is due to the fact that the push-forward of a uniform probability measure through a linear map gives a uniform push-forward measure with constant density.  Thus, the normalizing constant in the statistical Bayesian approach is the push-forward of the prior.
However, we expect to see significant difference between the statistical and consistent posterior densities for a nonlinear map.

We evaluate the model using 100,000 samples from the prior probability distribution.
For the consistent Bayesian approach, we use a standard KDE to approximate the push-forward of the prior
and use rejection sampling to select a subset of these samples for the posterior.
To make the comparison as close as possible, we evaluate the statistical Bayesian likelihood for each of these samples and estimate the normalizing constant which is simply the integral of the likelihood with respect to the prior.  We then use rejection sampling to select a subset of these samples for the statistical Bayesian posterior.  We obtained comparable results using a Markov chain Monte Carlo methods to generate samples, but these results are not reported here.
We then use a KDE to estimate the push-forward of each posterior using the samples accepted from the rejection sampling.

\begin{figure}[ht]
 \begin{center}
 \includegraphics[width=0.4\textwidth]{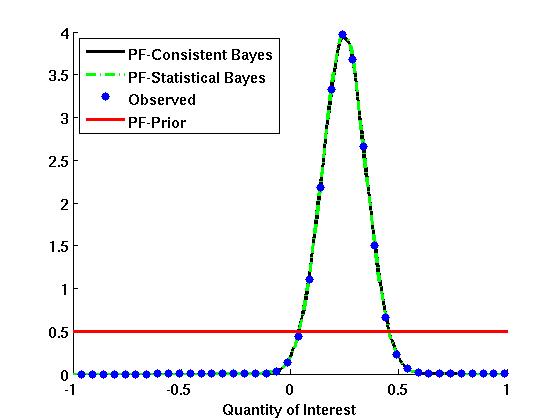}
\includegraphics[width=0.4\textwidth]{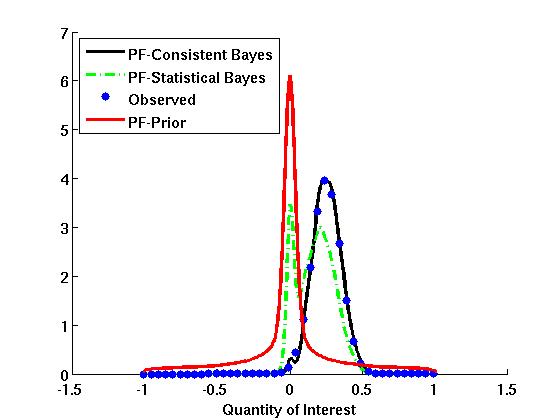}
 \end{center}
 \caption{The push-forward of the prior and the posteriors using statistical Bayesian inference and consistent Bayesian inference for a linear model (right) and for a nonlinear model (left).}
 \label{fig:ex1_postcomp}
 \end{figure}

In Figure~\ref{fig:ex1_postcomp}, we plot the push-forward of the prior, the push-forward of the consistent Bayesian posterior, and the push-forward of the statistical Bayesian posterior for $p=1$ (left) and $p=5$ (right) with $\bar{q} = 0.25$ and $\sigma = 0.1$.
For the linear map, the push-forward of both posteriors match the observed density.
However, we see that for the nonlinear map, the push-forward of the consistent Bayesian posterior is in reasonable agreement with the observed density while
the push-forward of the statistical Bayesian posterior is clearly a combination of the observed density and the push-forward of the prior.
This illustrates the fact that these two approaches solve different problems, give different posteriors and make different predictions.


\subsection{The Measure-theoretic Approach}\label{subsec:compmeas}

In this section, we briefly summarize the measure-theoretic approach and results first developed and analyzed in \cite{BE1, BE2, BE3} using notation in this work to make comparisons more obvious.  
The stochastic inverse problem formulated within the measure-theoretic framework is to determine a probability measure $P_\pspace$ on $(\pspace, \pborel)$ that is consistent with a given $\obsmeas$ on $(\dspace, \dborel)$ in the same sense as the consistency requirement of Eq.~\eqref{eq:invdefn}, i.e., the push-forward measure of $P_\pspace$ should agree with $\obsmeas$.
We do not refer to $P_\pspace$ as a posterior measure since it is not defined using the typical Bayesian formulation.


When $Q$ is at least piecewise smooth, $m<n$, and the Jacobian of $Q$ has full rank a.e. in $\Lambda$, then for a fixed datum $q\in\dspace$,  $Q^{-1}(q)$  exists as a (piecewise defined) $(n-m)$-dimensional manifold, which is referred to as a {\em generalized contour} in \cite{BE1,BE2,BE3}. 
Furthermore, a (possibly piecewise defined) $m$-dimensional manifold, called a transverse parameterization, exists in $\pspace$ that intersects each generalized contour exactly once.
This results in a type of equivalence class representation of points in $\pspace$ in terms of the generalized contours to which they belong.
Transverse parameterizations are not unique, and a fixed transverse parameterization is simply a particular explicit representation of this equivalence class defined by specifying a representative element from each equivalence class.
Constructing these manifolds is not necessary, but it is useful for understanding the structure of solutions to the stochastic inverse problem and for comparing $P_\Lambda$ to the posterior measure obtained in this work.

Following \cite{BE1,BE2,BE3}, let $\mathcal{L}$ define the space of equivalence classes in $\pspace$ where a particular $\ell\in\mathcal{L}$ corresponds to an entire generalized contour in $\pspace$. 
The map $Q$ is then a bijection between $\mathcal{L}$ and $\dspace$. 
This implies that $(\mathcal{L}, \mathcal{B}_{\mathcal{L}})$ is a measurable space where the $\sigma$-algebra is induced by $\dborel$ and $Q$. 
Furthermore, since $Q$ is a bijection between $\mathcal{L}$ and $\dspace$, we have that any probability measure $P_\dspace$ given on $(\dspace, \dborel)$ defines a unique probability measure on $(\mathcal{L}, \mathcal{B}_{\mathcal{L}})$ that is consistent with $P_\dspace$. 
Thus, given $\obsmeas$, there exists a {\em unique} probability measure $P_\pspace$ defined on $(\pspace, \mathcal{C}_\pspace)$ that is consistent with the model and $\obsmeas$.
However, the goal of solving the stochastic inverse problem is to determine $P_\pspace$ defined on $(\pspace, \pborel)$ not a measure on $(\pspace, \mathcal{C}_\pspace)$. 
To this end, \cite{BE1,BE2,BE3} uses the disintegration theorem and an ansatz to obtain such a measure. 

In the absence of any additional information, the so-called standard ansatz, first introduced in \cite{BE3}, uses the disintegration of the volume measure $\pmeas$ to define a family of uniform probability measures on the (uncountable) family of generalized contours that partition $\pspace$ into the equivalence classes.
Assuming the measures are described in terms of densities (with respect to various volume measures or their disintegrations), the end result is that computing probabilities using the density $\pi_\pspace$ describing $P_\pspace$ can be written in terms of iterated integrals using the marginal density, $\pi_{\mathcal{L}}(\text{proj}_{\mathcal{L}}(\lambda))$, uniquely defined by $\obsdens$, and the family of conditional densities on each generalized contour denoted by $\set{\pi_{\ell}(\lambda)}_{\ell\in\mathcal{L}}$ given by the ansatz.
Here, $\text{proj}_{\mathcal{L}}$ denotes a projection map from points in $\pspace$ to their equivalence class given in $\mathcal{L}$.  

We now remark on some specific similarities and differences between the solution to the consistent Bayesian approach and the measure $P_\Lambda$ described above.
The consistency requirement implies that the disintegration of $\postmeas$ on any transverse parameterization representing $\mathcal{L}$ must be equal to $\pi_{\mathcal{L}}(\text{proj}_{\mathcal{L}}(\lambda))$. 
Thus, if $\priormeas$ is defined so that its disintegration leads to the same family of conditional probabilities on the generalized contours as used by the ansatz in the measure-theoretic approach, then $P_\Lambda$ and $\postmeas$ will be identical due to the uniqueness of disintegrations \cite{BE3, Dellacherie_Meyer}.
In the case of the standard ansatz being applied, this equality can only be true when $\priormeas$ is a uniform probability measure with respect to the volume measure $\pmeas$.
In other words, when the $\priormeas$ is defined simply by a scaling of $\pmeas$, then the disintegration of $\priormeas$ gives a family of uniform probability measures along each generalized contour. 
When the $\priormeas$ is not uniform, it is unclear how to formulate a non-standard ansatz without explicitly using $\priormeas$ and its disintegration to construct a $P_\Lambda$ that matches identically $\postmeas$ on $(\pspace, \pborel)$. 
Moreover, the techniques used to construct an approximation to $P_\Lambda$ were based on approximations of {\em events} and the probabilities of these events with the equivalence class representation $\mathcal{L}$ allowing for all computations and disintegrations to take place explicitly in $\pspace$. 
This is in stark contrast to the techniques used here to construct a {\em density} approximation to $\postmeas$ based on sampling directly from the measure without first explicitly constructing its approximation and uses the standard form of the disintegration requiring evaluations of samples in both $\pspace$ and $\dspace$. 


In summary, the consistency requirement implies that the marginals of $\postmeas$ and $P_\Lambda$ on any transverse parameterization are identical, and any discrepancies between these measures on $(\pspace, \pborel)$ are due to discrepancies in $\priormeas$ and the ansatz.
Consequently, if $Q$ is a bijection between $\pspace$ and $\dspace$, then necessarily $\postmeas$ and $P_\Lambda$ always agree. 
In general, $Q$ is not a bijection so we expect to see differences between $\postmeas$ and $P_\Lambda$ but the push-forward of both
will agree with the observed measure.
Moreover, the numerical approaches used to approximate the density of $\postmeas$ and the measure $P_\Lambda$ are quite different due to the different assumptions and formulations.
In particular, the measure $P_\Lambda$ must be approximated on $\pspace$, which is often high-dimensional, while many operations on $\postdens$ can be done using only $r(Q(\lambda))$ which only involves approximations in $\dspace$ which is often low-dimensional.

\section{Conclusion}\label{sec:conclusions}
In this paper, we have introduced the {\em consistent} Bayesian formulation for 
stochastic inverse problems and proven that the posterior probability measure 
given by this approach is consistent with the model and the data.  
We also derived a simple expression for the posterior probability density and 
gave several numerical examples that use this expression to solve the stochastic inverse problem.
Comparisons were made with existing approaches, namely the statistical Bayesian approach 
and the measure theoretic approach.  
We gave an example where the consistent and statistical
give very different results and we discussed the scenarios where the consistent Bayesian and measure-theoretic approaches would give identical results.
We note that our approach does not require additional computations if more data
becomes available, i.e., if the observed measure/density changes, and is therefore
very amenable for Bayesian optimal experimental design.
Future work will incorporate more computationally efficient procedures to approximate the
push-forward of the prior probability density and to adaptively generate samples from the posterior.

\section{Acknowledgments}

T.~Butler's work is supported by the National Science Foundation (DMS-1228206). J.D.~Jakeman's work was partially supported by DARPA EQUIPS.

\bibliographystyle{plain}
\bibliography{references}

\begin{thebibliography}{10}

\bibitem{Bernardo1994}
{\sc J.~M. Bernardo and F.~M. Adrian}, {\em {Bayesian Theory}}, Wiley, 1994.

\bibitem{Billingsly}
{\sc P.~Billingsley}, {\em Probability and Measure}, John Wiley \& Sons, Inc.,
  1995.

\bibitem{BE1}
{\sc J.~Breidt, T.~Butler, and D.~Estep}, {\em A computational measure
  theoretic approach to inverse sensitivity problems {I}: Basic method and
  analysis}, SIAM J. Numer. Analysis, 49 (2012), pp.~1836--1859.

\bibitem{Bui-ThanhGhattasMartinEtAl13}
{\sc T.~Bui-Thanh, O.~Ghattas, J.~Martin, and G.~Stadler}, {\em A computational
  framework for infinite-dimensional {B}ayesian inverse problems {P}art {I}:
  {T}he linearized case, with application to global seismic inversion}, SIAM
  Journal on Scientific Computing, 35 (2013), pp.~A2494--A2523.

\bibitem{bungartz04}
{\sc H.-J. Bungartz and M.~Griebel}, {\em Sparse grids}, Acta Numerica, 13
  (2004), pp.~147--269.

\bibitem{BE2}
{\sc T.~Butler, D.~Estep, and J.~Sandelin}, {\em A measure-theoretic
  computational method for inverse sensitivity problems {II}: A posterior error
  analysis}, SIAM J. Numer. Analysis, 50 (2012), pp.~22--45.

\bibitem{BE3}
{\sc T.~Butler, D.~Estep, S.~Tavener, C.~Dawson, and J.~J. Westerink}, {\em A
  measure-theoretic computational method for inverse sensitivity problems iii:
  Multiple quantities of interest}, SIAM/ASA Journal on Uncertainty
  Quantification, 2 (2014), pp.~174--202.

\bibitem{Dellacherie_Meyer}
{\sc C.~Dellacherie and P.~Meyer}, {\em Probabilities and Potential},
  North-Holland Publishing Co., Amsterdam, 1978.

\bibitem{Devroye85}
{\sc L.~Devroye and L.~Gy{\''o}rfi}, {\em {Nonparametric Density Estimation:
  The L$_1$ View}}, Wiley, New York, 1985.

\bibitem{Ferenczi1997}
{\sc S.~Ferenczi}, {\em Measure-theoretic complexity of ergodic systems},
  Israel Journal of Mathematics, 100 (1997), pp.~189--207.

\bibitem{ganis2008stochastic}
{\sc B.~Ganis, H.~Klie, M.~F. Wheeler, T.~Wildey, I.~Yotov, and D.~Zhang}, {\em
  Stochastic collocation and mixed finite elements for flow in porous media},
  Computer methods in applied mechanics and engineering, 197 (2008),
  pp.~3547--3559.

\bibitem{Gelman2013}
{\sc A.~Gelman, J.~B. Carlin, H.~S. Stern, D.~B. Dunson, A.~Vehtari, and D.~B.
  Rubin}, {\em {Bayesian Data Analysis, Third Edition}}, Chapman and Hall/CRC,
  2013.

\bibitem{gerstner03}
{\sc T.~Gerstner and M.~Griebel}, {\em {Dimension-adaptive tensor-product
  quadrature}}, {Computing}, {71} ({2003}), pp.~{65--87}.

\bibitem{GhanemRedhorse}
{\sc R.~Ghanem and J.~Red-Horse}, {\em Propagation of probabilistic uncertainty
  in complex physical systems using a stochastic finite element approach},
  Physica D: Nonlinear Phenomena, 133 (1999), pp.~137--144.

\bibitem{GhanemSpanos}
{\sc R.~Ghanem and P.~Spanos}, {\em Stochastic {F}inite {E}lements: {A Spectral
  Approach}}, Springer Verlag, New York, 2002.

\bibitem{ChoosingMetrics_Gibbs}
{\sc A.~L. Gibbs and F.~E. Su}, {\em On choosing and bounding probability
  metrics}, INTERNAT. STATIST. REV.,  (2002), pp.~419--435.

\bibitem{hegland2003}
{\sc M.~Hegland}, {\em Adaptive sparse grids}, in Proc. of 10th Computational
  Techniques and Applications Conference CTAC-2001, K.~Burrage and R.~B. Sidje,
  eds., vol.~44, 2003, pp.~C335--C353.

\bibitem{Jakeman2013790}
{\sc J.~D. Jakeman, A.~Narayan, and D.~Xiu}, {\em Minimal multi-element
  stochastic collocation for uncertainty quantification of discontinuous
  functions}, Journal of Computational Physics, 242 (2013), pp.~790 -- 808.

\bibitem{Jaynes1998}
{\sc E.~T. Jaynes}, {\em {Probability Theory: The Logic of Science}}, 1998.

\bibitem{Kaipio_S_ACIP_2007}
{\sc J.~Kaipio and E.~Somersalo}, {\em Statistical inverse problems:
  Discretization, model reduction and inverse crimes}, Journal of Computational
  and Applied Mathematics, 198 (2007), pp.~493 -- 504.
\newblock Special Issue: Applied Computational Inverse Problems.

\bibitem{MaitreGhanem1}
{\sc O.~Le~Ma\^{i}tre, R.~Ghanem, O.~Knio, and H.~Najm}, {\em Uncertainty
  propagation using wiener-haar expansions}, J. Comput. Phys., 197(1) (2004),
  pp.~28--57.

\bibitem{Ma:2009:AHS:1514432.1514547}
{\sc X.~Ma and N.~Zabaras}, {\em An adaptive hierarchical sparse grid
  collocation algorithm for the solution of stochastic differential equations},
  J. Comput. Phys., 228 (2009), pp.~3084--3113.

\bibitem{MarzoukNajmRahn}
{\sc Y.~Marzouk, H.~Najm, and L.~Rahn}, {\em Stochastic spectral methods for
  efficient {B}ayesian solution of inverse problems}, J. Comput. Physics, 224
  (2007), pp.~560--586.

\bibitem{Narayan_GX_SISC_2013}
{\sc A.~Narayan, C.~Gittelson, and D.~Xiu}, {\em A stochastic collocation
  algorithm with multifidelity models}, SIAM Journal on Scientific Computing,
  36 (2014), pp.~A495--A521.

\bibitem{Ng_E_AIAA_2012}
{\sc L.~W.-T. Ng and M.~Eldred}, {\em Multifidelity uncertainty quantification
  using nonintrusive polynomial chaos and stochastic collocation}, in
  Proceedings of the 14th AIAA Non-Deterministic Approaches Conference, number
  AIAA-2012-1852, Honolulu, HI, vol.~43, 2012.

\bibitem{rasmussen2006}
{\sc C.~E. Rasmussen and C.~K.~I. Williams}, {\em Gaussian processes for
  machine learning},  (2006).

\bibitem{Robert2001}
{\sc C.~P. Robert}, {\em {The Bayesian Choice - A Decision Theoretic Motivation
  (second ed.)}}, Springer, 2001.

\bibitem{Royset_W_EJOR_2015}
{\sc J.~O. Royset and R.~J.-B. Wets}, {\em Fusion of hard and soft information
  in nonparametric density estimation}, European Journal of Operational
  Research, 247 (2015), pp.~532 -- 547.

\bibitem{Sargsyan_NG_IJCK_2015}
{\sc K.~Sargsyan, H.~N. Najm, and R.~Ghanem}, {\em On the statistical
  calibration of physical models}, International Journal of Chemical Kinetics,
  47 (2015), pp.~246--276.

\bibitem{Schwab2006100}
{\sc C.~Schwab and R.~A. Todor}, {\em Karhunen-{L}o\'eve approximation of
  random fields by generalized fast multipole methods}, Journal of
  Computational Physics, 217 (2006), pp.~100 -- 122.
\newblock Uncertainty Quantification in Simulation Science.

\bibitem{Stuart_AN_2010}
{\sc A.~M. Stuart}, {\em Inverse problems: A {B}ayesian perspective}, Acta
  Numerica, 19 (2010), pp.~451--559.

\bibitem{Tarantola}
{\sc A.~Tarantola}, {\em Inverse Problem Theory and Methods for Model Parameter
  Estimation}, SIAM, 2005.

\bibitem{Terrel_S_JSTOR_1992}
{\sc G.~R. Terrell and D.~W. Scott}, {\em Variable kernel density estimation},
  The Annals of Statistics, 20 (1992), pp.~1236--1265.

\bibitem{WanKarniadakis}
{\sc X.~Wan and G.~Karniadakis}, {\em {Beyond {W}iener-{A}skey Expansions:
  Handling Arbitrary PDFs}}, Journal of Scientific Computing, 27 (2006),
  pp.~455--464.

\bibitem{wand1994multivariate}
{\sc M.~Wand and M.~Jones}, {\em Multivariate plug-in bandwidth selection},
  Computational Statistics, 9 (1994), pp.~97--116.

\bibitem{wheeler2011multiscale}
{\sc M.~F. Wheeler, T.~Wildey, and I.~Yotov}, {\em A multiscale preconditioner
  for stochastic mortar mixed finite elements}, Computer Methods in Applied
  Mechanics and Engineering, 200 (2011), pp.~1251--1262.

\bibitem{XiuKarniadakis}
{\sc D.~Xiu and G.~Karniadakis}, {\em The {W}iener-{A}skey polynomial chaos for
  stochastic differential equations}, SIAM J. Sci. Comput., 24 (2002),
  pp.~619--644.

\bibitem{zhang2004efficient}
{\sc D.~Zhang and Z.~Lu}, {\em An efficient, high-order perturbation approach
  for flow in random porous media via {K}arhunen-{L}o\'eve and polynomial
  expansions}, Journal of Computational Physics, 194 (2004), pp.~773--794.

\end{thebibliography}

\end{document}